\documentclass[a4paper,12pt]{amsart}

\usepackage[hyperindex]{hyperref}

\usepackage{t1enc}
\def\frak{\mathfrak}
\def\Bbb{\mathbb}
\def\Cal{\mathcal}

\def\F{\mathbb{F}}
\def\L{\mathbb{L}}
\def\R{\mathbb{R}}
\def\W{\mathbb{W}}
\def\X{\mathbb{X}}
\def\Y{\mathbb{Y}}
\def\Z{\mathbb{Z}}
\def\N{\mathbb{N}}

\def\cE{\mathcal{E}}
\def\cV{\mathcal{V}}
\def\cT{\mathcal{T}}

\def\cC{\mathcal{C}}
\def\cH{\mathcal{H}}
\def\cN{\mathcal{N}}
\def\cL{\mathcal{L}}
\def\cR{\mathcal{R}}
\def\cB{\mathcal{B}}

\def\de{\delta}

\def\la{\lambda}
\def\rh{\rho}
\def\si{\sigma}

\def\om{\omega}

\def\De{\Delta}

\def\Om{\Omega}
\def\Up{\Upsilon}
\def\na{\nabla}

\def\form#1{\mathbf{#1}}
\def\dform#1{\dot{\mathbf{#1}}}
\def\ddform#1{\ddot{\mathbf{#1}}}

\newcommand{\id}{\operatorname{id}}
\newcommand{\End}{\operatorname{End}}

\newcommand{\lpl}{
  \mbox{$
  \begin{picture}(12.7,8)(-.5,-1)
  \put(2,0.2){$+$}
  \put(6.2,2.8){\oval(8,8)[l]}
  \end{picture}$}}
\renewcommand{\vec}[1]{\mathbf{#1}}

\newcommand{\wh}{\widehat}
\newcommand{\cq}{{\Cal Q}}

\newcommand{\Proj}{\operatorname{Proj}}

\newcommand{\modDe}{\fl}
\newcommand{\modD}{\fD}
\newcommand{\modBox}  {\mbox{$                                    
\begin{picture}(9,8)(1.6,0.15)
\put(1,0.2){\mbox{$ \boxcircle \hspace{-7.8pt} /$}}
\end{picture}$}} 

\newcommand{\newc}{\newcommand}
\newcommand{\aU}{\mbox{\boldmath{$ U$}}}

\newtheorem{theorem}{Theorem}[section]
\newtheorem{lemma}[theorem]{Lemma}
\newtheorem{proposition}[theorem]{Proposition}
\newtheorem{corollary}[theorem]{Corollary}
\theoremstyle{remark}
\newtheorem{remark}[theorem]{\rm\bf Remark}
\newtheorem*{remark*}{\rm\bf Remark}

\newcommand{\ce}{{\Cal E}}

\usepackage{amssymb,stmaryrd}
\usepackage{amscd}

\newcommand{\nd}{\nabla}

\newcommand{\Rho}{P}

\newcommand{\nn}[1]{(\ref{#1})}

\newcommand{\D}{\mbox{\boldmath{$ D$}}}
\newcommand{\tU}{\tilde{U}}

\newcommand{\bg}{\mbox{\boldmath{$ g$}}}

\newcommand{\J}{J}

\newcommand{\NX}{{\mbox{\boldmath$ \nabla$}}_{\!\miniX}}
\newcommand{\aM}{\tilde{M}}
\newcommand{\sX}{\mbox{\scriptsize\boldmath{$X$}}}
\newcommand{\cce}{\tilde{\ce}}
\newcommand{\miniX}{\mbox{\boldmath{$\scriptstyle{X}$}}}
\newcommand{\ct}{{\Cal T}}
\newcommand{\act}{\mbox{\boldmath{${\mathcal{T}}$}}}
\newcommand{\bsi}{\mbox{\boldmath{$\si$}}}
\newcommand{\aD}{\mbox{\boldmath{$D$}}}

\newcommand{\fl}{\mbox{$
\begin{picture}(9,8)(1.6,0.15)
\put(1,0.2){\mbox{$ \Delta \hspace{-7.8pt} /$}}
\end{picture}$}}
\newcommand{\afl}{\mbox{$
\begin{picture}(9,8)(1.6,0.15)
\put(1,0.2){\mbox{$ {\mbox{\boldmath$\Delta$}} \hspace{-7.8pt} /$}}
\end{picture}$}}
\newcommand{\fD}{\mbox{$
\begin{picture}(9,8)(1.6,0.15)
\put(1,0.2){\mbox{$ D \hspace{-7.8pt} /$}}
\end{picture}$}}
\newcommand{\afD}
{\mbox{$
\begin{picture}(9,8)(1.6,0.15)
\put(1,0.2){\mbox{$ \D \hspace{-7.8pt} /$}}
\end{picture}$}}

\newcommand{\boldnabla}{\mbox{\boldmath$ \nabla$}}
\newcommand{\nda}{\boldnabla}
\newcommand{\aX}{\mbox{\boldmath{$ X$}}}
\newcommand{\aI}{\mbox{\boldmath{$ I$}}}
\newcommand{\af}{\mbox{\boldmath{$ f$}} }
\newcommand{\aDelta}{\mbox{\boldmath{$ \Delta$}}}

\newc{\aR}{\mbox{\boldmath{$ R$}}}
\newc{\aS}{\mbox{\boldmath{$ S$}}}
\newc{\aDeR}{\mbox{\boldmath{$ U$}}_B{}^P{}_C{}^Q}
\newc{\aDe}{\mbox{\boldmath$ \Delta$}}
\newc{\aNd}{\mbox{\boldmath$ \nabla$}}

\newc{\aK}{\mbox{\boldmath{$ K$}}}
\newc{\aL}{\mbox{\boldmath{$ L$}}}

\newcommand{\h}{\mbox{\boldmath{$ h$}}}

\def\sideremark#1{\ifvmode\leavevmode\fi\vadjust{\vbox to0pt{\vss
 \hbox to 0pt{\hskip\hsize\hskip1em
 \vbox{\hsize3cm\tiny\raggedright\pretolerance10000
 \noindent #1\hfill}\hss}\vbox to8pt{\vfil}\vss}}}%
                        
\def\idx#1{{\em #1\/}}

\author{A. Rod Gover and Josef \v Silhan}
\email{gover@math.auckland.ac.nz} \title{Conformal
operators on forms and detour complexes on Einstein manifolds}

\begin{document}

\begin{abstract}
For even dimensional conformal manifolds several new conformally
invariant objects were found recently: invariant differential
complexes related to, but distinct from, the de Rham complex (these
are elliptic in the case of Riemannian signature); the cohomology
spaces of these; conformally stable form spaces that we may view as
spaces of conformal harmonics; operators that generalise Branson's
Q-curvature; global pairings between differential form bundles that
descend to cohomology pairings.  Here we show that these operators,
spaces, and the theory underlying them, simplify significantly on
conformally Einstein manifolds. We give explicit formulae for all the
operators concerned. The null spaces for these, the conformal
harmonics, and the cohomology spaces are expressed explicitly in terms
of direct sums of subspaces of eigenspaces of the form Laplacian.  For
the case of non-Ricci flat spaces this applies in all signatures and
without topological restrictions.  In the case of Riemannian signature
and compact manifolds, this leads to new results on the global
invariant pairings, including for the integral of Q-curvature against
the null space of the dimensional order conformal Laplacian of Graham
et al..
\end{abstract}

\address{ARG: Department of Mathematics\\
  The University of Auckland\\
  Private Bag 92019\\
  Auckland 1\\
  New Zealand} \email{gover@math.auckland.ac.nz}
\address{JS: Eduard \v{C}ech Center \\ 
Department of Algebra and geometry \\
Masaryk University \\
Jan\'a\v{c}kovo n\'am. 2a \\
602 00, Brno\\
Czech Republic} \email{silhan@math.muni.cz}

\maketitle

\pagestyle{myheadings}
\markboth{Gover \& \v Silhan}{Conformal geometry of forms on Einstein manifolds}

\section{Introduction}

Differential forms provide a fundamental domain for the study of
smooth manifolds. In Riemannian geometry the de Rham complex, its
associated Hodge theory, and distinguished forms representing
characteristic classes are among the most basic and important tools
(e.g.\ \cite{CS74,deR}). In physics the study of forms is partly
motivated by Maxwell theory and its generalisations. Operators on
differential forms feature strongly in string and brane theories. In
both mathematics and physics Einstein manifolds have a central position
\cite{Besse} and thus they give an important class of special
structures for the study of geometric objects.

Among the differential operators that are natural for
pseudo--Rieman\-nian structures only a select class are conformally
invariant.  Conformal invariance is a subtle property which reflects
an independence of point dependent scale. This symmetry is manifest in
the equations of massless particles. It is linked to CR geometry (hence
complex analysis) through the Fefferman metric \cite{F}; the natural
equations on the Fefferman space are conformally invariant.  This
symmetry also underpins the conformal approach to Riemannian geometry.
For example, it is essentially exploited in the Yamabe problem (see
\cite{Schoen} and references therein) of prescribing the scalar
curvature. Recently there has been a focus on variations of this idea,
including the conformal prescription of Branson's Q-curvature
\cite{tomsharp,CQY,malchiodi}. These problems use the conformal Laplacian
on functions (or densities) and its higher order analogues due to
Paneitz, Graham et al.\ \cite{GJMS}.

The use of conformal operators on forms provides a setting where, on
the one hand, there is potential to formally generalise such theories,
but which, on the other hand, should yield access to rather different
geometric data. An immediate difficulty is that forms are more
difficult to work with than functions and so, while there was much
early work in this direction (e.g. \cite{tbms,tomgp}), this did not yield a
clear picture. In dimension 4, and inspired by constructions from
twistor theory, some rather interesting directions and applications to
global geometry were pioneered in the work of Eastwood and Singer
\cite{EastSin1,EastSin2}. Links between this result and the tractor
calculus of \cite{BEGo,CapGotrans,luminy} were established in
\cite{BrGoMax}. On the other hand in \cite{CapGoamb,GoPetLap} it is
shown that the conformal tractor connection may be recovered as a
suitable linearisation of the ambient metric of Fefferman and Graham
\cite{FGast} (and see also \cite{FGrNew}). Exploiting both developments a rather complete theory of
conformal operators on forms was derived in the joint works
\cite{BrGodeRham,BrGocomm} of the first author with Branson. The main point of
that article was not simply to construct conformal operators on
differential forms, but rather, to expose and develop the discovery of
preferred versions of such operators and the rather elegant picture
that these yield: one may immediately construct, on even dimensional
conformal manifolds, a host of new global conformal invariants. Some of
these generalise, in a natural way, the integral of Q-curvature.

For most of these new operators explicit formulae are not available. For
any particular operator a formula may be obtained algorithmically
via tractor calculus and the theory developed in
\cite{GoPetLap,GoPetobst}. However the resulting operators, when
presented in the usual way, are given by extremely complicated
formulae. It turns out there are striking simplifications when these
operators are studied on conformally Einstein manifolds. The purpose
of this article is to expose this, via a comprehensive but concise treatment,
and use the results to study, in the Einstein setting, the related  global 
conformal invariants and spaces.

To describe the content in more detail we first review  the
relevant results from \cite{BrGodeRham} and \cite{BrGocomm}.  On
conformal manifolds of even dimension $n\geq 4$ there is a family of
formally self-adjoint conformally invariant differential complexes:
\begin{equation}\label{dcx}
\ce^0\stackrel{d}{\to}\cdots\stackrel{d}{\to}\ce^{k-1}\stackrel{d}{\to}
\ce^k\stackrel{L_k}{\to}\ce_k\stackrel{\delta}{\to}\ce_{k-1}
\stackrel{\delta}{\to}
\cdots\stackrel{\delta}{\to}\ce_0~.
\end{equation}
 Here, for each $k\in\{0,1,\cdots ,n/2-1\}$, $\ce^k$ denotes the space of $k$-forms, $\ce_k$
denotes an appropriate density twisting of that space, $d$ is the exterior
derivative and $\delta$ its formal adjoint. 
An interesting feature of these complexes is that the operators  
$L_k$ have
the structure of a composition 
$$
L_k=\delta Q_{k+1} d
$$ where $Q_{k+1}$ is from a family of
differential operators, parametrised by $k=-1,\cdots ,n/2-1$, and
which, as operators on closed forms, generalise Branson's Q-curvature;
in particular under conformal rescaling of the metric $g\mapsto
\widehat{g}=e^{2\om}g$ ($\om\in C^\infty$) these have the conformal
transformation formula
$$
\hat Q_k\,u=Q_ku+L_k(\om u)\ \ 
\mbox{for }u\mbox{ a {\em closed} }k\mbox{-form},
$$ $Q_0 1$ is the Q-curvature and $L_0$ is the dimension order GJMS
operator of \cite{GJMS}.  On closed forms these {\em Q-operators} have
the form $Q_{k+1}= (d\delta)^{n/2-k-1}+lower~order~terms$, so in the
case of Riemannian signature the complexes \nn{dcx} are
elliptic. Writing $H^k_L$ for the (conformally invariant) cohomology
at $k$, for the complex \nn{dcx}, it follows that on compact
Riemannian manifolds $H^k_L$ is finite.

The composition $G_k := \delta Q_k$ is a conformal gauge companion operator
for $L_k$ and also for the exterior derivative $d$. What this means is
that the systems $(L_k,G_k)$ and $(d,G_k)$ are, in a suitable sense,
conformally invariant and, in the case of Riemannian signature, are
graded injectively elliptic. For example the null space
${\mathcal{H}}_G^k$ of $(d,G_k)$ is conformally stable and, as pointed
out in \cite{BrGodeRham}, is a candidate for a space of ``conformal
harmonics''. Some perspective on these objects is given by the sequence
\begin{equation}\label{coco}
0\to H^{k-1}\to H^{k-1}_L\stackrel{d}{\to}\mathcal{H}^k_G \to H^k, \quad \quad k\in\{1,2,\cdots ,n/2\}
\end{equation}
where $H^k$ indicates the usual de Rham cohomology. The map
$H^{k-1}\to H^{k-1}_L$ is the obvious inclusion, since $L$ factors
through the exterior derivative. The map $\mathcal{H}^k_G \to H^k$ is
that which simply takes solutions of $(d,G_k)$ to their cohomology
class. It is immediate from the definitions of the spaces that the
sequence is exact, but it is an open question whether in general the
final map $\mathcal{H}^k_G \to H^k $ is surjective. When it is we say that
the space is $(k-1)$-regular \cite{BrGodeRham}.

We present here a study of all of these spaces and operators
specialised to the setting of an Einstein structure. By exploiting
some recent developments we obtain a treatment which, surprisingly,
obtains most of the results in a uniform way in all signatures and
without assuming the manifold is compact. As indicated above the
motivation is manifold.  The cohomology spaces, and related structures,
mentioned above are clearly fundamental to conformal geometry. An
important problem is to discover what data they capture.
On the other hand there is the opportunity to shed light on Einstein
structures which form an important class of geometries which remain
rather mysterious; for example there are very few non-existence
results for compact Riemannian Einstein spaces, while the construction of
examples is primarily through K\"ahler geometry.  
The idea that this might be a rewarding approach is suggested by the
intimate relationship between conformal geometry and Einstein
structures.  Conformal structures admit a natural conformally
invariant connection on a prolonged structure: this is the Cartan
connection of \cite{Cartan}, or equivalently the induced structure is
the conformal standard tractor connection that was already
mentioned. An Einstein structure is equivalent to a suitably generic
parallel section of this tractor bundle and so is, in this sense, a
type of symmetry of conformal structure.

For case of the conformal Laplacian type operators this last point was
exploited heavily in \cite{GoEinst} where two of the three main
results are as follows: on Einstein manifolds the GJMS operators of
\cite{GJMS} factor into compositions of operators each of which is of
the form of a constant potential Helmholtz Laplacian; the Q-curvature
is constant and (up to a universal constant) simply a power of the
scalar curvature.  (See also \cite{FGrNew} where similar results are
obtained using different techniques.) Here we develop the analogous
theory for operators on differential forms. In fact we do much more.
A first step is that we obtain factorisations of the key operators
which generalise those from \cite{GoEinst} (to our knowledge such
factorisations are new even for the conformally flat Einstein
setting).  On the other hand in \cite{GoSiDec,GoSiSrni} we show that if the factors
$P_i:\cV\to \cV$ (for some vector space $\cV$), in a composition
$P:=P_0P_1 \cdots P_\ell$ of mutually commuting operators, are
suitably ``relatively invertible'' then the general inhomogeneous
problem $Pu=f$ decomposes into an equivalent system $P_iu_i=f$,
$i=0,\cdots , \ell$. This is used extensively in the current work to
reduce, on non-Ricci flat Einstein manifolds, the generally high order
conformal operators to equivalent lower order systems. The outcome is
that in any signature (and without any assumption of compactness) on
non-Ricci flat Einstein manifolds we can describe the spaces
$\cN(L_k)$ (the null space of $L_k$), and $\cH^k_G$ explicitly as a
direct sum of $ \cH^k_{\si} := \cN(d) \cap \cN(\de)$ and the (possibly
trivial) ``eigenspaces'',
\begin{eqnarray*}
  &\overline{\cH}^k_{\si,\la} := \{ f \in \cE^k \mid d\de f = \la f\},
   \quad
   \widetilde{\cH}^k_{\si,\la} := \{ f \in \cE^k \mid \de d f = \la f\} ,
\end{eqnarray*}
for various explicitly known $\lambda \in \mathbb{R}$.
(Here $\si$ denotes the Einstein scale in the conformal class,
see Section \ref{null}.)  See in
particular Proposition \ref{nullL} and \nn{decomp}, and note
that for $\lambda\neq 0$ the displayed spaces give the $\cR(d)$ (range
of $d$) and $\cR(\de)$ parts of the form Laplacian ``eigenspace'' $ \{
f \in \cE^k \mid (d\de+\de d ) f = \la f\}$. 
We also come to a simple
decomposition for $H^k_L$ (see expression \nn{Hkl}) and other
conformal spaces from \cite{BrGodeRham}.

Stronger results are available in the compact Riemannian setting and
these are summarised in Theorem \ref{crnrf}. Observe, in particular,
that this shows that all compact Riemannian Einstein even manifolds
are $k$-regular for $k=0,1\cdots,n/2-1$, and that, in this setting,
$\cH^k_\si$ agrees with the usual space of harmonics for the form
Laplacian.  From the $k$-regularity it follows that the global
conformally invariant pairings on $\cH^k_G$, as defined in
\cite{BrGocomm}, descend to a conformal quadratic form on de Rham
cohomology. See Theorem \ref{Rpairs}, and also Proposition \ref{Qdes}
which shows that, in the Einstein case, the pairing is given by a
power of the scalar curvature; in fact by a formula which generalises
the formula from \cite{GoEinst} for Q-curvature on Einstein
manifolds. Also in Theorem \ref{Rpairs} we show that the conformal
pairing, via Q-operators, of closed forms against forms in $\cN(L_k)$
descends to a closed form pairing. See the Remark following Theorem
\ref{Rpairs}, which emphasises that this also gives a result for the
usual Q-curvature.  Some of the results for compact Riemannian
manifolds could be obtained by using, at the outset, the complete spectral
resolution of the form Laplacian. However doing this conceals the fact
that, for the most part, the same results are available even when we
do not have access to diagonalisations of the basic operators.

The development is as follows. Section \ref{tractorsect} summarises
some basic conformal geometry, tractor results, and identities to be
used.  In Section \ref{ops} we construct Laplacian operators on
weighted tractor bundles. This is in the spirit of \cite{GoEinst}, but
there is an algebraic adjustment to the basic operators.  Using these
Laplacian power operators, in section \ref{detour} we derive formulae
for the key operators, $L_k$, $Q_k$, and so forth, in the Einstein
setting. The main result is Theorem \ref{mainthm}.  In fact in
contrast to the construction in \cite{BrGodeRham} (which heavily uses
the Fefferman-Graham ambient metric), these operators are developed
and defined directly using invariant tractor operators. That we
recover the operators from \cite{BrGodeRham} is the the main subject
of Section \ref{asec}.  In each case the operators are given in terms
of compositions of commuting operators. This enables, in Section
\ref{null}, the use of the tools from \cite{GoSiDec} as recounted in
Theorem \ref{fundthm}.

The first author would like to thank the Royal Society of New Zealand
for support via Marsden Grant no. 06-UOA-029.  The second author was
supported from Basic Research Center no.\ LC505 (Eduard \v{C}ech
Center for Algebra and Geometry) of the Ministry of Education of Czech
Republic. We are appreciative of the careful reading by the referee;
this exposed a number of typographical errors in the original manuscript.

\section{Background: Einstein metrics and conformal geometry} 

We first sketch here notation and background for general conformal
structures and their tractor calculus following
\cite{CapGoamb,GoPetLap}. The latter is then used to describe
operators that we will need acting on tractor forms and some key
identities are developed. Some parts of the treatment are specialised to
Einstein manifolds.

\subsection{Conformal geometry and tractor calculus}\label{tractorsect}

 Let $M$ be a smoo\-th manifold of
dimension $n\geq 3$. Recall that a {\em conformal structure\/} of
signature $(p,q)$ on $M$ is a smooth ray subbundle $\cq\subset
S^2T^*M$ whose fibre over $x$ consists of conformally related
signature-$(p,q)$ metrics at the point $x$. Sections of $\cq$ are
metrics $g$ on $M$. The principal bundle $\pi:\cq\to M$ has structure group
$\Bbb R_+$, and each representation ${\Bbb R}_+ \ni x\mapsto
x^{-w/2}\in {\rm End}(\Bbb R)$ induces a natural line bundle on $
(M,[g])$ that we term the conformal density bundle $E[w]$. We shall
write $ \ce[w]$ for the space of sections of this bundle. Here and
throughout, sections, tensors, and functions are always smooth.  When
no confusion is likely to arise, we will use the same notation for a
bundle and its section space.

We write $\bg$ for the {\em conformal metric}, that is the
tautological section of $S^2T^*M\otimes E[2]$ determined by the
conformal structure. This will be used to identify $TM$ with
$T^*M[2]$.  For many calculations we will use abstract indices in an
obvious way.  Given a choice of metric $ g$ from the conformal class,
we write $ \nabla$ for the corresponding Levi-Civita connection. With
these conventions the Laplacian $ \Delta$ is given by
$\Delta=\bg^{ab}\nd_a\nd_b= \nd^b\nd_b\,$.  Note $E[w]$ is trivialised
by a choice of metric $g$ from the conformal class, and we write $\nd$
for the connection arising from this trivialisation.  It follows
immediately that (the coupled) $ \nd_a$ preserves the conformal
metric.

Since the Levi-Civita connection is torsion-free, the (Riemannian)
curvature 
$R_{ab}{}^c{}_d$ is given by $ [\nd_a,\nd_b]v^c=R_{ab}{}^c{}_dv^d $ where 
$[\cdot,\cdot]$ indicates the commutator bracket.  The Riemannian
curvature can be decomposed into the totally trace-free Weyl curvature
$C_{abcd}$ and a remaining part described by the symmetric {\em
Schouten tensor} $\Rho_{ab}$, according to $
R_{abcd}=C_{abcd}+2\bg_{c[a}\Rho_{b]d}+2\bg_{d[b}\Rho_{a]c}, $ where
$[\cdots]$ indicates antisymmetrisation over the enclosed indices.
We shall write $J := P^a{}_a$. The {\em
Cotton tensor} is defined by
$$
A_{abc}:=2\nabla_{[b}\Rho_{c]a} .
$$
Under a {\em conformal transformation} we replace a choice of metric $
g$ by the metric $ \hat{g}=e^{2\om} g$, where $\omega$ is a smooth
function. Explicit formulae for the corresponding transformation of
the Levi-Civita connection and its curvatures are given in e.g.\ 
\cite{GoPetLap}. We recall that, in particular, the Weyl curvature is
conformally invariant $\widehat{C}_{abcd}=C_{abcd}$.

We next define the standard tractor bundle over $(M,[g])$.
It is a vector bundle of rank $n+2$ defined, for each $g\in[g]$,
by  $[\ce^A]_g=\ce[1]\oplus\ce_a[1]\oplus\ce[-1]$. 
If $\wh g=e^{2\Up}g$, we identify  
 $(\alpha,\mu_a,\tau)\in[\ce^A]_g$ with
$(\wh\alpha,\wh\mu_a,\wh\tau)\in[\ce^A]_{\wh g}$
by the transformation
\begin{equation}\label{transf-tractor}
 \begin{pmatrix}
 \wh\alpha\\ \wh\mu_a\\ \wh\tau
 \end{pmatrix}=
 \begin{pmatrix}
 1 & 0& 0\\
 \Up_a&\delta_a{}^b&0\\
- \tfrac{1}{2}\Up_c\Up^c &-\Up^b& 1
 \end{pmatrix} 
 \begin{pmatrix}
 \alpha\\ \mu_b\\ \tau
 \end{pmatrix} ,
\end{equation}
where $\Up_a:=\nd_a \Up$.  It is straightforward to verify that these
identifications are consistent upon changing to a third metric from
the conformal class, and so taking the quotient by this equivalence
relation defines the {\em standard tractor bundle} $\ct$, or $\ce^A$
in an abstract index notation, over the conformal manifold.
(Alternatively the standard tractor bundle may be constructed as a
canonical quotient of a certain 2-jet bundle or as an associated
bundle to the normal conformal Cartan bundle \cite{luminy}.) On a
conformal structure of signature $(p,q)$, the bundle $\ce^A$ admits an
invariant metric $ h_{AB}$ of signature $(p+1,q+1)$ and an invariant
connection, which we shall also denote by $\nabla_a$, preserving
$h_{AB}$.  In a conformal scale $g$, these are given by
\begin{equation}\label{basictrf}
 h_{AB}=\begin{pmatrix}
 0 & 0& 1\\
 0&\bg_{ab}&0\\
1 & 0 & 0
 \end{pmatrix}
\text{ and }
\nabla_a\begin{pmatrix}
 \alpha\\ \mu_b\\ \tau
 \end{pmatrix}
 =
\begin{pmatrix}
 \nabla_a \alpha-\mu_a \\
 \nabla_a \mu_b+ \bg_{ab} \tau +\Rho_{ab}\alpha \\
 \nabla_a \tau - \Rho_{ab}\mu^b  \end{pmatrix}. 
\end{equation}
It is readily verified that both of these are conformally well-defined,
i.e., independent of the choice of a metric $g\in [g]$.  Note that
$h_{AB}$ defines a section of $\ce_{AB}=\ce_A\otimes\ce_B$, where
$\ce_A$ is the dual bundle of $\ce^A$. Hence we may use $h_{AB}$ and
its inverse $h^{AB}$ to raise or lower indices of $\ce_A$, $\ce^A$ and
their tensor products.

In computations, it is often useful to introduce 
the `projectors' from $\ce^A$ to
the components $\ce[1]$, $\ce_a[1]$ and $\ce[-1]$ which are determined
by a choice of scale.
They are respectively denoted by $X_A\in\ce_A[1]$, 
$Z_{Aa}\in\ce_{Aa}[1]$ and $Y_A\in\ce_A[-1]$, where
 $\ce_{Aa}[w]=\ce_A\otimes\ce_a\otimes\ce[w]$, etc.
 Using the metrics $h_{AB}$ and $\bg_{ab}$ to raise indices,
we define $X^A, Z^{Aa}, Y^A$. Then we
immediately see that 
$$
Y_AX^A=1,\ \ Z_{Ab}Z^A{}_c=\bg_{bc} ,
$$
and that all other quadratic combinations that contract the tractor
index vanish. 
In \eqref{transf-tractor} note that  
$\wh{\alpha}=\alpha$ and hence $X^A$ is conformally invariant. 

Given a choice of conformal scale, the {\em tractor-$D$ operator} 
$D_A\colon\ce_{B \cdots E}[w] \to \cE_{AB\cdots E}[w-1]$
is defined by 
\begin{equation}\label{tractorD}
D_A V:=(n+2w-2)w Y_A V+ (n+2w-2)Z_{Aa}\nabla^a V -X_A\Box V, 
\end{equation} 
where $\Box V :=\Delta V+w \J V$.  This also turns out to be
conformally invariant as can be checked directly using the formulae
above (or alternatively there are conformally invariant constructions
of $D$, see e.g.\ \cite{Gosrni}).

The curvature $ \Omega$ of the tractor connection 
is defined by 
$$
[\nd_a,\nd_b] V^C= \Omega_{ab}{}^C{}_EV^E 
$$
for $ V^C\in\ce^C$.  Using
\eqref{basictrf} and the formulae for the Riemannian curvature yields
\begin{equation}\label{tractcurv}
\Omega_{abCE}= Z_C{}^cZ_E{}^e C_{abce}-2X_{[C}Z_{E]}{}^e A_{eab}
\end{equation}

We will also need a conformally invariant curvature quantity defined as
follows (cf.\ \cite{Gosrni,Goadv})
\begin{equation}\label{Wdef}
W_{BC}{}^E{}_F:=
\frac{3}{n-2}D^AX_{[A} \Omega_{BC]}{}^E{}_F ,
\end{equation}
where $\Omega_{BC}{}^E{}_F:= Z_B{}^bZ_C{}^c \Om_{bc}{}^E{}_F$.
In a choice of conformal scale, 
 $W_{ABCE}$ is given by
\begin{equation}\label{Wform}
\begin{array}{l}
(n-4)\left( Z_A{}^aZ_B{}^bZ_C{}^cZ_E{}^e C_{abce}
-2 Z_A{}^aZ_B{}^bX_{[C}Z_{E]}{}^e A_{eab}\right. \\ 
\left.-2 X_{[A}Z_{B]}{}^b Z_C{}^cZ_E{}^e A_{bce} \right)
+ 4 X_{[A}Z_{B]}{}^b X_{[C} Z_{E]}{}^e B_{eb},
\end{array}
\end{equation}
where 
$$
B_{ab}:=\nabla^c
A_{acb}+\Rho^{dc}C_{dacb}.
$$ is known as the {\em Bach tensor}. From the formula \nn{Wform} it
is clear that $W_{ABCD}$ has Weyl tensor type symmetries.

We will work with conformally Einstein manifolds. That is, conformal
structures with an Einstein metric in the conformal class. This is the
same as the existence of a non-vanishing section $\si \in
\mathcal{E}[1]$ satisfying $\left[ \na_{(a}\na_{b)_0} + P_{(ab)_0}
\right] \si =0$, where $(\ldots)_0$ indicates the trace-free 
symmetric part over the enclosed indices. 
Equivalently (see e.g.\ \cite{BEGo,GoNur}) there is a standard
tractor $I_A$ that is parallel with respect to the normal tractor
connection $\na$ and such that $\si:=X_AI^A$ is non-vanishing. It
follows that $I_A := \frac{1}{n}D_A \si = Y_A \si +Z_A^a \na_a \si
-\frac{1}{n} X_A (\De+J) \si$, for some section $\si\in \ce[1]$, and
so $X^AI_A=\si$ is non-vanishing. If we compute in the scale $\si$,
then for example $W_{ABCD} = (n-4) Z_A^a Z_B^b Z_C^c Z_D^d C_{abcd}$.

\subsection{Tractor forms}
Following \cite{BrGodeRham} we
write $\cE^k[w]$ for the space of sections of 
$(\Lambda^k T^*M) \otimes E[w]$ (and $\ce^k= \ce^k[0]$).  
Further we put $\cE_k[w] := \cE^k[w+2k-n]$. The space of closed 
$k$-forms shall be denoted by $\cC^k \subseteq \cE^k$.

In order to be
explicit and efficient in calculations involving bundles of possibly
high rank it is necessary to employ abstract index notation
as follows.  
In the usual abstract index conventions one would write
$\ce_{[ab\cdots c]}$ (where there are implicitly $k$-indices skewed
over) for the space $\ce^k$. To simplify subsequent expressions we
 use the following conventions. Firstly
indices labelled with sequential superscripts which are
at the same level (i.e.\ all contravariant or all
covariant) will indicate a completely skew set of indices.
Formally we set $a^1 \cdots a^k = [a^1 \cdots a^k]$ and so, for example,
$\ce_{a^1 \cdots a^k}$ is an alternative notation for $\ce^k$
while $\ce_{a^1 \cdots a^{k-1}}$ and $\ce_{a^2 \cdots a^k}$ both denote
$\ce^{k-1}$. Next, following \cite{GoSiKil} we
abbreviate this notation via multi-indices: We will use the forms
indices
$$
\begin{aligned}
 \vec{a}^k &:=a^1 \cdots a^k =[a^1 \cdots a^k], \quad k \geq 0,\\
\dot{\vec{a}}^k &:= a^2 \cdots a^k=[a^2 \cdots a^k], \quad k \geq 1.
\end{aligned}
$$
If $k=1$ then $\dform{a}^k$ simply means the index is absent.
The corresponding notations will be used for tractor indices so
e.g. the bundle of tractor $k$--forms $\ce_{[A^1\cdots A^k]}$ will be
denoted by $\ce_{A^1\cdots A^k}$ or $\mathcal{E}_{\vec{A}^k}$.

The structure of $\mathcal{E}_{\vec{A}^k}$ is
\begin{equation} \label{comp_series_form}
  \mathcal{E}_{[A^1 \cdots A^k]} = \mathcal{E}_{\vec{A}^k} \simeq
  \mathcal{E}^{k-1}[k] \lpl \left( \mathcal{E}^k[k] \oplus
  \mathcal{E}^{k-2}[k-2] \right) \lpl \mathcal{E}^{k-1}[k-2];
\end{equation}
in a choice of scale the semidirect sums $\lpl$ may be replaced by
direct sums and otherwise they indicate the composition series
structure arising from the tensor powers of \nn{transf-tractor}.

 In 
a choice of metric $g$ from the conformal class, the
projectors (or splitting operators) $X,Y,Z$ for $\mathcal{E}_A$
determine corresponding projectors $\X,\Y,\Z,\W$ for
$\mathcal{E}_{\vec{A}^{k+1}}$, $k \geq 1$ 
These execute the  splitting of this space into four components and are given 
as follows.
\begin{center}
\renewcommand{\arraystretch}{1.3}
\begin{tabular}{c@{\;=\;}l@{\;=\;}l@{\;=\;}l@{\ $\in$\ }l}
$\Y^k$ & $\Y_{A^0A^1 \cdots A^k}^{\quad a^1 \cdots\, a^k}$ &
  $\Y_{A^0\vec{A}^k}^{\quad \vec{a}^k}$ & $Y_{A^0}^{}Z_{A^1}^{a^1} \cdots Z_{A^k}^{a^k}$ &
  $\mathcal{E}_{\vec{A}^{k+1}}^{\vec{a}^k}[-k-1]$ \\
$\Z^k$ & $\Z_{A^1 \cdots A^k}^{\, a^1 \cdots\, a^k}$ &
  $\Z_{\vec{A}^k}^{\,\vec{a}^k}$ & $Z_{A^1}^{\,a^1} \cdots Z_{A^k}^{\,a^k}$ &
  $\mathcal{E}_{\vec{A}^k}^{\vec{a}^k}[-k]$ \\
$\W^k$ & $\W_{A'A^0A^1 \cdots A^k}^{\quad\,\ \ a^1 \cdots\, a^k}$ &
  $\W_{A'A^0\vec{A}^k}^{\quad\,\ \ \vec{a}^k}$ &
  $X_{[A'}^{}Y_{A^0}^{}Z_{A^1}^{\,a^1} \cdots Z_{A^k]}^{\,a^k}$ &
  $\mathcal{E}_{\vec{A}^{k+2}}^{\vec{a}^k}[-k]$ \\
$\X^k$ & $\X_{A^0A^1 \cdots A^k}^{\quad a^1 \cdots\, a^k}$ &
  $\X_{A^0\vec{A}^k}^{\quad \vec{a}^k}$ & $X_{A^0}^{}Z_{A^1}^{\,a^1} \cdots Z_{\,A^k}^{a^k}$ &  
  $\mathcal{E}_{\vec{A}^{k+1}}^{\vec{a}^k}[-k+1]$
\end{tabular}
\end{center}
where $k \geq 0$. The superscript $k$ in $\Y^k$, $\Z^k$, $\W^k$
and $\X^k$ shows the corresponding tensor valence. (This is
slightly different than in \cite{BrGodeRham}, where $k$ is
the relevant tractor valence.) Note that $Y=\Y^0$, $Z=\Z^1$ and $X=\X^0$ and
$\W^0 = X_{[A'}Y_{A^0]}$.
From \nn{basictrf} we immediately see $\na_p Y_A = Z_A^a P_{pa}$, 
$\na_p Z_A^a = -\delta_p^a Y_A - P_p^a X_A$ and $\na_p X_A = Z_{Ap}$.
From this we obtain 
the formulae (cf.\ \cite{GoSiKil}) 
\begin{equation} \label{na}
\begin{split}
  \na_p \Y_{A^0\form{A}^k}^{\quad \form{a}^k} &=
      P_{pa_0} \Z_{A^0\form{A}^k}^{\,a^0\,\form{a}^k}
    + k P_p^{\ a^1} \W_{A^0\form{A}^k}^{\quad \dot{\form{a}}^k} \\
  \na_p \Z_{A^0\form{A}^k}^{\,a^0\,\form{a}^k} &=
    - (k+1) \de_p^{a^0} \Y_{A^0\form{A}^k}^{\quad \form{a}^k}
    - (k+1) P_p^{\ a^0}  \X_{A^0\form{A}^k}^{\quad \form{a}^k} \\
  \na_p \W_{A^0\form{A}^k}^{\quad\, \dot{\form{a}}^k} &=
    - \bg_{pa^1} \Y_{A^0\form{A}^k}^{\quad \form{a}^k}
    + P_{pa^1} \X_{A^0\form{A}^k}^{\,\ \ a^1\!\dot{\form{a}}^k} \\
  \na_p \X_{A^0\form{A}^k}^{\quad \form{a}^k} &=
      \bg_{pa^0} \Z_{A^0\form{A}^k}^{\,a^0\,\form{a}^k}
    - k \de_p^{a^1} \W_{A^0\form{A}^k}^{\quad \dot{\form{a}}^k},
\end{split}
\end{equation}
which determine the tractor connection on form tractors in a conformal
scale. Similarly, one can compute the Laplacian $\De$ applied to the
tractors $\X$, $\Y$, $\Z$ and $\W$.  As an operator on form tractors
we have the opportunity to modify $\De$ by adding some amount of
$W\sharp\sharp$, where $\sharp $ denotes the natural tensorial
action of sections in $\End(\ce^A)$. Analogously, we shall use
$C \sharp\sharp$ to modify the Laplacian on forms; here $\sharp$ 
denotes the natural tensorial action of sections in $\End(\ce^a)$.
It turns out
(cf.\ \cite{BrGodeRham}) that it will be convenient for us to  use the
operator
$$ \modDe = \begin{cases}
      \De + \frac{1}{n-4} W\sharp\sharp & n \not= 4 \\
      \De & n=4.
   \end{cases} $$
 (Note $\De = \na^a \na_a$.)
Since the Laplacian is of the second order, it is  convenient 
to consider e.g.\ 
$\modDe \Y_{\form{A}}^{\,\dform{a}} \tau_{\dform{a}}$ where 
$\tau_{\dform{a}} \in \cE_{\dform{a}}[w]$. It will be sufficient for our 
purpose
to calculate this only in an Einstein scale. 
For example, using \nn{na} and then that $P_{ab}=\bg_{ab}J/n$, we have 
\begin{align*}
\na^p \na_p \Y_{\form{A}}^{\,\dform{a}} \tau_{\dform{a}} 
   =& \na^p \bigl[ 
    P_{pa^1} \Z_{\form{A}}^{\form{a}}
    + (k-1) P_p^{\ a^2} \W_{\form{A}}^{\ddot{\form{a}}}
    + \Y_{\form{A}}^{\,\dform{a}} \na_p  \bigr]\tau_{\dform{a}} \\
   =&-\Y_{\form{A}}^{\,\dform{a}}
     \bigl[ \bigl( \de d + d\de + (1\!-\!\frac{2(k\!-\!1)(n\!-\!k\!+\!1)}{n}) J
             +C \sharp\sharp \bigr)\tau \bigr]_{\dform{a}}   \\
    &+\frac{2}{nk} \Z_{\form{A}}^{\,\form{a}} (Jd\tau)_{\form{a}}
     -\frac{2(k\!-\!1)}{n} \W_{\form{A}}^{\,\ddform{a}} (J\de\tau)_{\ddform{a}} 
     -\frac{n\!-\!2k\!+\!2}{n^2} \X_{\form{A}}^{\,\dform{a}}
       J^2\tau_{\dform{a}}, 
\end{align*}
where, as usual, $\form{A}=\form{A}^k$ and $\form{a}=\form{a}^k$.
Summarising, one can compute that in
an Einstein scale we obtain 
\begin{equation} \label{dde-W}
\begin{split}
  - \modDe \Y_{\form{A}}^{\,\dform{a}} \tau_{\dform{a}} &= 
     \Y_{\form{A}}^{\,\dform{a}}
     \bigl[ \bigl( \de d + d\de + (1\!-\!\frac{2(k\!-\!1)(n\!-\!k\!+\!1)}{n}) J
            \bigr)\tau \bigr]_{\dform{a}}   \\
    &-\frac{2}{nk} \Z_{\form{A}}^{\,\form{a}} (Jd\tau)_{\form{a}}
     +\frac{2(k\!-\!1)}{n} \W_{\form{A}}^{\,\ddform{a}} (J\de\tau)_{\ddform{a}} 
     +\frac{n\!-\!2k\!+\!2}{n^2} \X_{\form{A}}^{\,\dform{a}}
       J^2\tau_{\dform{a}} \\
  - \modDe \Z_{\form{A}}^{\,\form{a}} \mu_{\form{a}} &= 
     -2k \Y_{\form{A}}^{\,\dform{a}} (\de\mu)_{\dform{a}} 
     +\Z_{\form{A}}^{\,\form{a}} \bigl[ \bigl(
       \de d + d\de - \frac{2k(n\!-\!k\!-\!1)}{n} J
       \bigr)\mu \bigr]_{\form{a}} \\
    &-\frac{2k}{n} \X_{\form{A}}^{\,\dform{a}} (J\de\mu)_{\dform{a}} \\
  - \modDe \W_{\form{A}}^{\,\ddform{a}} \nu_{\ddform{a}} &= 
       \frac{2}{k\!-\!1} \Y_{\form{A}}^{\,\dform{a}} (d\nu)_{\dform{a}} 
    +\W_{\form{A}}^{\,\ddform{a}} \bigl[ \bigl(
      \de d + d\de - \frac{2(k\!-\!3)(n\!-\!k\!+\!2)}{n}  J
       \bigr) \nu \bigr]_{\ddform{a}} \\
    &-\frac{2}{n(k\!-\!1)} \X_{\form{A}}^{\,\dform{a}} (d\nu)_{\dform{a}} \\
  - \modDe \X_{\form{A}}^{\,\dform{a}} \rh_{\dform{a}} &= 
     (n\!-\!2k\!+\!2) \Y_{\form{A}}^{\,\dform{a}} \rh_{\dform{a}} 
     -2(k\!-\!1) \W_{\form{A}}^{\,\ddform{a}} (\de\rh)_{\ddform{a}} \\ 
    &-\frac{2}{k} \Z_{\form{A}}^{\,\form{a}} (d\rh)_{\form{a}} 
     +\X_{\form{A}}^{\,\dform{a}} \bigl[ \bigl(
       \de d + d\de + (1\!-\!\frac{2(k\!-\!1)(n\!-\!k\!+\!1)}{n}) J
       \bigr) \rh \bigr]_{\dform{a}}.
\end{split}
\end{equation}
if either $n=4$, $k=1$ or $n \not= 4$, cf.\ \cite[(1.50)]{Sthesis}.
Here $\tau_{\dform{a}} \in \cE_{\dform{a}}[w]$, 
$\mu_{\form{a}} \in \cE_{\form{a}}[w]$, $\nu_{\ddform{a}} \in \cE_{\ddform{a}}[w]$
and $\rh_{\dform{a}} \in \cE_{\dform{a}}[w]$ where $\form{a} = \form{a}^k$,
$k \geq 1$ and $w$ is any conformal weight.

\subsection{Useful identities}\label{ids}

Here we first introduce and discuss some identities that hold on a
general conformal manifold.

Recall that sequentially labelled indices are assumed to be skew over,
e.g.\ $A^1A^2 = [A^1A^2]$.
The operator 
\begin{equation} \label{doubleD}
  D_{A^1A^2} = -2( w\W_{A^1A^2} + \X_{A^1A^2}^{\quad a} \na_a)
\end{equation}
was introduced in \cite{Goadv}. Also recall the definition of the tractor
$W_{B^1B^2C^1C^2}$ in \nn{Wdef}. By replacing $\De$ by $\modDe$ in
\nn{tractorD} we obtain the conformally invariant operator
$$ \modD_A = \begin{cases}
     D_A - \frac{1}{n-4} X_A W\sharp\sharp & n \not= 4 \\
     D_A & n=4. 
   \end{cases} 
$$
The case $n \not= 4$ of this operator was introduced in 
\cite{BrGodeRham}. 
In contrast to $D$, and surprisingly, the commutator of $\modD$ is algebraic (cf.\ 
the commutator of $D$ in \cite{GoPetLap}) for $n \not= 4$:
\begin{gather} \label{[modD,modD]}
  [\modD_A,\modD_B] = \frac{n+2w-2}{n-4} \bigl[
  (n+2w-4) W_{AB} \sharp  - (D_{AB} W) \sharp\sharp \bigr] f,\ n \not= 4. 
\end{gather}
This can be checked by a direct computation, or alternatively by a
rather simple calculation using the ambient metric and its links to
tractors as in Section \ref{asec}.  For $n=4$ we have
$[\modD_A,\modD_B] = [D_A,D_B]$, see \cite{GoPetLap} for the latter.
Note one can moreover show that for $n \not= 4$, the operator $D_A - y
X_A W\sharp\sharp$, $y \in \R$ has algebraic commutator only for the
value $y = \frac{1}{n-4}$.

\begin{proposition} \label{doubleDW}
$D_{A^1A^2}$ and $W_{B^1B^2C^1C^2}$ have the following properties: \\
(i) $D_{[A^1A^2} W_{B^1B^2]C^1C^2} =0$  \\
(ii) $D_{A^1}{}^P W_{PA^2B^1B^2} = -W_{A^1A^2B^1B^2}$.
\end{proposition}

\begin{proof}
We shall use the form indices $\form{A} = \form{A}^2$, $\form{B} = \form{B}^2$ 
and $\form{C} = \form{C}^2$ throughout the proof.
Both identities can be verified by the direct computation.
To simplify the computation note that alternatively \cite{Sthesis} we have
$$ W_\form{BC} = \bigl[ (n-4) \Z_{\form{B}}^\form{b} 
   - 2 \X_{\form{B}}^{b^2} \na^{\,b^1} \bigr] \Om_\form{bC}
   \in \cE_\form{BC}[-2]. $$

(i) Using the relations $\W_{[\form{A}}^{} \X_{\form{B}]}^{\,b} =
\X_{[\form{A}}^{\,a} \X_{\form{B}]}^{\,b} = 
\X_{[\form{A}}^{\,b} \W_{\form{B}]}^{} =0$ 
(which follow from $X_{[A} X_{B]}=0$) we obtain
\begin{align*} 
  D_{[\form{A}} W_{\form{B}]\form{C}} = 
  & -2(n-4) \W_{[\form{A}}^{} \Z_{\form{B}]}^{\form{b}} \Om_\form{bC} 
    -2(n-4) \X_{[\form{A}}^{\ a} \Y_{\form{B}]}^{\,b} \Om_{ab\form{C}} \\
  &  +(n-4) \X_{[\form{A}}^{\,a} \Z_{\form{B}]}^{\form{b}} \na_a \Om_{\form{bC}}
    -2 \X_{[\form{A}}^{\ a} \Z_{\form{B}]}^{\form{b}} 
     \bg_{ab^1} \na^p \Om_{pb^2\form{C}}.
\end{align*}
Now clearly the first two terms on the right hand side add up to $0$ and the 
remaining ones both vanish.

(ii) Clearly $\W_{A^1}{}^P \Z_{PA^2}^{\,a^1a^2} =  
\X_{A^1}{}^{Pa^1} \X_P^{}{}_{A^2}^{a^2} =0$. Thus
\begin{align*} 
  D_{A^1}{}^P W_{PA^2\form{B}} = -2\Bigl\{
  & 4(n-4) \W_{A^1}{}^P \X_P^{}{}_{A^2}^{a^2} \na^q\Om_{qa^2\form{B}} \\
  & + (n-4) \X_{A^1}^{}{}^{Pp} \bigl[ 
    -2\Y_P^{}{}_{A^2}^{a^2} \Om_{pa^2\form{B}} 
    + \Z_{PA^2}^{\,a^1a^2} \na_p \Om_{\form{aB}} \bigr] \\
  & + \X_{A^1}^{}{}^{Pp} \bigl[
    -2 \Z_{PA^2}^{\,a^1a^2} \bg_{pa^1} + 2 \W_{PA^2} \de_p^{a^2} \bigr]
    \na^q \Om_{qa^2\form{B}}
  \Bigr\}.
\end{align*}
Now using $\W_{A^1}{}^P \X_P^{}{}_{A^2}^{a^2} = 
\frac{1}{4} \X_{\form{A}}^{\,a^2}$, $\X_{A^1}^{}{}^{Pp} Z_{PA^2}^{\,a^1a^2} = 
\frac{1}{2} \X_{\form{A}}^{\,a^2} \bg^{pa^1}$ and
$\X_{A^1}^{}{}^{Pp} \Y_P^{}{}_{A^2}^{a^2} = 
-\frac{1}{4} \Z_{\form{A}}^{pa^2} - \frac{1}{4} \W_\form{A} \bg^{pa^2}$
the proposition follows.
\end{proof}

\begin{lemma} \label{IdoubleDW}
(i) Let $I^A \in \cE^A$ be a parallel tractor. Then
$I^P D_P{}^Q W_{QA^2B^1B^2} =0$ and 
$I^P D_{P[A^0} W_{A^1A^2]B^1B^2} =0$. 

(ii) Let $I^A, \bar{I}^A \in \cE^A$ be a two parallel tractors. Then
$I^{A^1} \bar{I}^{A^2} D_{A^1A^2} W_{B^1B^2C^1C^2} =0$. 
\end{lemma}

\begin{proof}
Recall any parallel tractor $I^A$ satisfies $I^{A^1} W_{A^1A^2B^1B^2}=0$
\cite{GoNur}. Thus the first relation of (i) follows by applying $I^{A^1}$
to Proposition \ref{doubleDW} (ii) and the second relation of (i) 
by applying $I^{A^1}$
to Proposition \ref{doubleDW} (i).
Similarly, (ii) follows by applying $I^{A^1} \bar{I}^{A^2}$
to Proposition \ref{doubleDW} (i).
\end{proof}

\section{Einstein manifolds: conformal Laplacian operators on tractors}
\label{ops}

We assume that the structure $(M,[g])$ is conformally Einstein, and write 
$\si \in \cE[1]$ for some  Einstein scale from
the conformal class. Then $I^A := \frac{1}{n} D^A \si$ is parallel and
$X^AI_A=\si$ is nonvanishing.

The operator $\Box = \De + wJ$ acting on tractor bundles of the weight
$w$ is invariant only if $n+2w-2=0$. On the other hand the scale $\si$
(or equivalently $I^A$), yields the operator
\begin{equation} \label{modBox}
  \modBox_\si := I^A \modD_A: \cE_{B \cdots E}[w] \longrightarrow 
  \cE_{B \cdots E}[w-1]
\end{equation}
which is well defined for any $w$, cf.\ \cite{GoEinst}. 
Thus we can consider
the composition $(\modBox_\si)^p$, $p \in \N$ and we set
$(\modBox_\si)^0 := \id$. These operators generally depend on the choice
of the Einstein scale but one has the following modification of
\cite[Theorem 3.1]{GoEinst}.

\begin{theorem} \label{indep}
Let $\si,\bar{\si}$ be two Einstein scales in the conformal
class and consider the operators 
$$ \frac{1}{\si^p} (\modBox_\si)^p, 
   \frac{1}{\bar{\si}^p} (\modBox_{\bar{\si}})^p:
   \cE_{B \cdots E}[w] \longrightarrow \cE_{B \cdots E}[w-2p],  $$
for $p\in \mathbb{Z}_{\geq 0}$.
If $w= p- n/2 $ then 
$\frac{1}{\si^p} (\modBox_\si)^p = 
\frac{1}{\bar{\si}^p} (\modBox_{\bar{\si}})^p$.
\end{theorem}

\begin{proof}
Assume $w= p- n/2 $ and denote the Einstein metric corresponding
to $\si$ by $g := \si^{-2}\bg$. Then $\frac{1}{\si^p} (\modBox_\si)^p$
is very similar to the operator $P_p^g$ given by \cite[(19)]{GoEinst}.
The difference is, beside the sign, that we have replaced $D_A$ in 
the definition of $P_p^g$ by $\modD_A$ in the definition of 
$\frac{1}{\si^p} (\modBox_\si)^p$. Thus our statement is analogous to 
\cite[Theorem 3.1]{GoEinst} for the operators $P_p^g$. We can follow
the proof of the latter theorem literally; the difference between
$D_A$ and $\modD_A$ appears only in the commutator in the last display
in the proof of \cite[Theorem 3.1]{GoEinst}. In our case, we need 
the commutator $[\modD_A,\modD_B]$ instead of that display and to 
finish the proof it remains to show that 
$I^A\bar{I}^B [\modD_A,\modD_B]$ vanishes on density valued tractor 
fields. Here $I^A = \frac{1}{n}D^A \si$, 
$\bar{I}^B = \frac{1}{n} D^B \bar{\si}$ are parallel 
tractors corresponding to the Einstein scales $\si$, $\bar{\si}$.
But this follows from \nn{[modD,modD]}, Lemma \ref{IdoubleDW} (ii) 
and the fact $I^{A^1}W_{A^1A^2B^1B^2}=0$.
\end{proof}

\section{Einstein
manifolds: Q-operators, gauge operators and detour complexes}
\label{detour}

The main aim here is to recover, in the Einstein setting, the
differential complexes \nn{dcx}, the Q-operators and the related
conformal spaces, as defined in \cite{BrGodeRham}. In that source the
Fefferman-Graham ambient metric is used to generate the operators
which form the ``building blocks'' of the theory. In contrast here, in
the conformally Einstein setting, all the operators are
``rediscovered'' directly using the tractor operators. However in
Section \ref{asec} we use the ambient metric to establish that we do
have exactly the specialisation of the operators and spaces from
\cite{BrGodeRham}; see the comment following expression \nn{LLk} and
Proposition \ref{ambLLk}.

Here we work on an Einstein manifold in an Einstein scale $\si \in \cE[1]$.
The first step is the conformally invariant differential splitting 
operator
\begin{eqnarray} \label{M}
\begin{split}
  &M_{\form{A}}^{\form{a}}: \cE_{\form{a}} \longrightarrow
     \cE_{\form{A}}[-k] \\
  &M_{\form{A}}^{\form{a}} f_{\form{a}} =
     \frac{n-2k}{k} \Z_{\form{A}}^{\,\form{a}} f_{\form{a}}
     + \X_{\form{A}}^{\,\dform{a}} (\de f)_{\dform{a}}.
\end{split}
\end{eqnarray}
where $\form{A} = \form{A}^k$ and $\form{a} = \form{a}^k$. 
This is a special case of the operator $\overline{M}$ from \cite{GoSiKil}
(up to the multiple $k$).

Let $I^A := \frac{1}{n}D^A \si$ be the Einstein tractor 
for $\si \in \cE[1]$.
Since
$I^A = Y^A \si - \frac{1}{n} X^A \si J$ in the scale $\si$, it follows from
\nn{modBox} that
\begin{equation} \label{EmodBox}
 \modBox_\si = \si (-\modDe - 2\frac{w}{n}(n+w-1)J): 
   \cE_{B \cdots E}[w] \longrightarrow \cE_{B \cdots E}[w-1] 
\end{equation}
in the scale $\si$.
For any differential operator $E: \cE^k \to \cE^k$ we have 
the compositions
\begin{equation} \label{Yam}
  P^p_k[E] := \prod_{i=1}^p
  \left( E + \frac{2i(n-2k-i+1)}{n}J \right)
\end{equation}
for $p \geq 1$. We set $P^p_k[E] := \id$ for $p \leq 0$.
Note that $P^p_k$ can be considered as a polynomial in $E$.

Next we define the operator
\begin{equation} \label{L}
  \L_k^p := \frac{1}{\si^p} (\modBox_\si)^p 
  M_{\form{A}}^{\form{a}}:
  \cE_\form{a}[w] \longrightarrow  \cE_\form{A} [w-k-2p]
\end{equation}
for $p \geq 0$ where $(\modBox_\si)^0:=\id$.
We put $\L_k := \L_k^p$ for $p = \frac{n+2w-2k}{2}$.
It follows from Theorem \ref{indep} that $\L_k$ is independent on 
the choice of the Einstein scale $\si$.

Now we are ready to state the main technical step of our construction.

\begin{theorem} \label{formula}
Let $f_{\form{a}} \in \cE_{\form{a}}$ where $\form{a} = \form{a}^k$
and $p \geq 1$.
Then computing in the Einstein scale $\si$ we obtain
\begin{align*}
  (\L_k^p f)_{\form{a}} =&
    -p(n-2k-2p) \Y_{\form{A}}^{\,\dform{a}} 
    [ \de P^{p-1}_k [d\de] f ]_{\dform{a}} \\
  &+\frac{1}{k} \Z_{\form{A}}^{\,\form{a}}
    [ (n-2k-2p) d\de P^{p-1}_k [d\de] f 
    + (n-2k) \de d P^{p-1}_{k+1} [\de d] f ]_{\form{a}} \\
  &+  \X_{\form{A}}^{\,\dform{a}}
    [\de (d\de + \frac{p(n-2k+2)}{n} J) P^{p-1}_k [d\de] f ]_{\dform{a}}.
\end{align*}
\end{theorem}

\begin{proof}
It is easy to show by a direct computation (using (\ref{na}) and (\ref{dde-W})) 
that the theorem holds for $p=1$. Now assume, by induction, that the theorem
holds for a fixed  $p \in \N$. To verify the theorem for $p+1$ we need
to compute
\begin{align*}
  (\L_k^{p+1} f)_{\form{a}} 
  &= \frac{1}{\si^{p+1}} \modBox_\si  (\modBox_\si)^p 
     M_{\form{A}}^{\form{a}} f_{\form{a}} \\
  &= \frac{1}{\si^{p+1}} 
     \si \bigl( -\modDe + 2 \frac{k+p}{n} (n-k-p-1)J \bigr)
     (\modBox_\si)^p M_{\form{A}}^{\form{a}} f_{\form{a}} \\
  &= \bigl( -\modDe + 2 \frac{k+p}{n} (n-k-p-1)J \bigr)
     (\L_k^p f)_{\form{a}}      
\end{align*}
where the multiple of $J$ follows from (\ref{EmodBox}) and 
from the conformal weight of
$(\modBox_\si)^p M_{\form{A}}^{\form{a}} f_{\form{a}} \in 
\cE_{\form{A}}[-k-p]$.
Since the theorem yields the formula for
$(\L_k^p f)_{\form{a}}$
we can continue in the computation using (\ref{na}) and (\ref{dde-W}).
The rest of the proof is to finish this computation. Let us show at least
the top slot. Using (\ref{na}), (\ref{dde-W}) and the formula for
$(\L_k^p f)_{\form{a}}$
we obtain 
\begin{align*}
  \Y_{\form{A}}^{\dform{a}} \Bigl[
  & - p (n-2k-2p) \de d \de P^{p-1}_k [d\de] f \\
  & - p (1-\frac{2(k-1)(n-k+1)}{n})(n-2k-2p) J\de P^{p-1}_k [d\de] f \\
  & - 2(n-2k-2p) \de d \de Y(p-1) f \\
  & +  (n-2k+2) \de (d\de+\frac{p(n-2k+2)}{n}J) P^{p-1}_k [d\de] f \\
  & - 2 \frac{k+p}{n} (n-k-p-1) p (n-2k-2p) J\de P^{p-1}_k [d\de] f 
      \Bigr]_{\dform{a}}.
\end{align*}
Summing up appropriate terms, we obtain that this is exactly
$$ -(p+1)(n-2k-2(p+1)) \Y_{\form{A}}^{\dform{a}} 
   [ \de P^p_k [d\de] f ]_{\dform{a}}. $$ 
The computation of remaining slots is similar and the theorem follows.
\end{proof}

Now let us assume that the dimension $n$ is even 
and $k \in
\{1,\ldots,\frac{n}{2}\}$. 
Then putting $p = \frac{n-2k}{2} \geq 0$
in Theorem \ref{formula}, we obtain
\begin{equation} \label{LLk}
  (\L_kf)_\form{a} =     
  \frac{2p}{k} \Z_{\form{A}}^{\,\form{a}}
  \bigl( \de d P^{p-1}_{k+1} [\de d] f \bigr)_\form{a}
  +\X_{\form{A}}^{\,\dform{a}}
  \bigl( \de P^p_k [d\de] f \bigr)_{\dform{a}}.
\end{equation}
for $f_\form{a} \in \cE_{\form{a}}$.
It follows from Proposition \ref{ambLLk} that
$\L_k$ acting on $\cE_\form{a}$ agrees with the operator $\L_k$ 
defined in \cite{BrGodeRham} up to a nonzero scalar multiple. 
Thus we have the following
results for the operators $G_k^\si$, $Q_k^\si$ and 
$L_k$ from \cite{BrGodeRham}.

\begin{theorem}\label{mainthm}
In an Einstein scale $\si$ we have, up to a nonzero scalar multiple,
the formula
\begin{equation}\label{Gk}
\de 
\prod_{i=1}^{\frac{n-2k}{2}} \left( d\de + \frac{2i(n-2k-i+1)}{n}J \right). 
\end{equation}
for the operator $G^\si_k$ of \cite{BrGodeRham}
$$
G^\si_k:\ce^k\to \ce_{k-1} .
$$
As an operator on closed $k$-forms, up to a nonzero scalar multiple, we have:
\begin{gather} \label{Q_k}
  Q^\si_k = 
 \prod_{i=1}^{\frac{n-2k}{2}} \left( d\de + \frac{2i(n-2k-i+1)}{n}J \right).
\end{gather}
\noindent
Note that at the $k=n/2$ extreme these are taken to mean $G_{n/2}^\si = \de$ and $Q_{n/2}^\si=\id$.
\end{theorem}

\begin{proof}
  Note that setting $p=\frac{n-2k}{2}$ in Theorem \ref{formula}, the
  coefficient of $\Y$ vanishes in the display and the left-hand-side
  of the display is the operator $\mathbb{L}_k$ of \cite{BrGodeRham}.
   The coefficient of $\X $ is thus
  $G_k^\si$ as in Theorem 4.5 of \cite{BrGodeRham}. This gives the
  formula presented.  
Here we have used the fact that the factor 
$d\de + \frac{p(n-2k+2)}{n} J$, with $p=\frac{n-2k}{2}$,  of the  coefficient of $\X$
appears as a composition factor in $P_k^p[d\de]$ (the factor with $i=p$ in \nn{Yam}).

Now by Theorem 2.8 of \cite{BrGodeRham} and its proof we have
  that $G_k^\si=\de Q^\si_k$ and $Q_k^\si$ is formally self-adjoint.
  Thus the formal adjoint $G_{k}^{\si,*}$ of $G_{k}^\si$ is $Q^\si_k
  d:\ce^{k-1}\to \ce_k$. On the other hand since $Q^\si_k$ is a differential
  operator it follows from this that $G_{k}^{\si,*}$ determines the given 
  formula for $Q^\si_k$ as an operator on closed forms.
\end{proof}

Observe that from \nn{Q_k} we have immediately the following useful
observation.
\begin{corollary} \label{QCC}
On Einstein manifolds and in an Einstein
scale $Q^\si_k: \mathcal{C}^k \to \mathcal{C}^k $.
\end{corollary}

The operator $G_k^\si$ acting on closed forms will be denoted by
$G_k: \cC^k \to \cE_{k-1}$. It follows from \nn{LLk} that 
$G_k$ is conformally invariant (recall it is defined as projection to 
the $\X$-slot of \nn{LLk}).

\begin{corollary} \label{L_k}
In an Einstein scale $\si$ the conformally
invariant detour operator $L_k:\ce^k\to \ce_k$ is given, up to a nonzero 
scalar multiple, by
$$
L_k = \de 
\prod_{i=1}^{\frac{n-2k-2}{2}} \left( d\de + \frac{2i(n-2k-i-1)}{n}J \right)
d 
$$
for $k<n/2$. Moreover we set $L_{n/2}=0$.
\end{corollary}

\begin{proof}
On a general manifold we have $L_k =\de Q_{k+1}^\si d$ 
from Theorem 2.8 of \cite{BrGodeRham}. Hence the statement
follows from the formula $Q_{k+1}^\si$ from \nn{Q_k}.
\end{proof}
\begin{remark} Observe that the operators $L_k$, and the fact that they 
have the form $\delta M d$, may be extracted directly from expression
\nn{LLk}.  Thus, in this conformally Einstein setting, we obtain
detour complexes as in \nn{dcx} from (\ref{LLk}). Of course these were
developed on general conformal manifolds in \cite{BrGodeRham} but the
construction in the conformally Einstein setting here is independent
of \cite{BrGodeRham}. Proposition \ref{ambLLk} is only used here to
verify that it is (the specialisation of) the same complex and also to
make the connection to the Q-operators from \cite{BrGodeRham}.
\end{remark}

We have constructed the operator $\L_k$ (thus also $G_k^\si$, 
$Q_k^\si$ and $L_k$)
for $k \in \{1,\ldots,\frac{n}{2}\}$ as we assumed the latter 
range in \nn{LLk}. However formulae \nn{Gk} for $G_k^\si$, 
\nn{Q_k} for $Q_k^\si$ and Corollary \ref{L_k} for $L_k$ make 
sense also for $k=0$ \cite{GoEinst};  the operators
$Q_0^\si$ and $L_0$ formally agree, up to a nonzero scalar multiple, with
the corresponding operators $Q^g$ and $\Box_{n/2}$ from \cite{GoEinst}.
Further we put $\L_0 := L_0$.
Here and below $g = \si^{-2} \bg$ is the Einstein metric.
So henceforth
we shall assume $k \in \{0,\ldots,\frac{n}{2}\}$.

\section{Decompositions of the conformal spaces}\label{null}

We work on a general (possibly noncompact) conformally Einstein even
dimensional manifold $(M,[g])$ of signature $(p,q)$.
As usual we write $\si$ to denote an Einstein scale.  If not stated
otherwise, we assume $k \in \{0,\ldots,\frac{n}{2}\}$ and we put
$\cE_{-1} := 0$.  The space $\cH^k_G = \cN(G_k: \cC^k \longrightarrow
\cE_{k-1})$, where $G_k := \de Q_k^\si$, is conformally invariant
\cite{BrGodeRham}, and we shall term it the space of {\em conformal
harmonics}.

We shall describe $\cH^k_G$ in more 
details. As mentioned in the Introduction, we will use the notation
\begin{eqnarray*}
  &\overline{\cH}^k_{\si,\la} := \{ f \in \cE^k \mid d\de f = \la f\},
   \quad
   \widetilde{\cH}^k_{\si,\la} := \{ f \in \cE^k \mid \de d f = \la f\} \\
  &\text{and} \quad \cH^k_{\si} := \cN(d) \cap \cN(\de)
\end{eqnarray*}
where $\la \in \R$.  Then $\cH^k_{\si} \subseteq
\overline{\cH}^k_{\si,0} + \widetilde{\cH}^k_{\si,0}$.  Note that if
$\la\neq 0 $ then $\overline{\cH}^k_{\si,\la}\subset \cR(d)$, and similarly 
$\widetilde{\cH}^k_{\si,\la}\subset \cR(\delta) $.

As well as $\cH^k_G$, we shall also study the null spaces of
the operators $G_k$ and $L_k$.  Our treatment relies on the following observation:

\begin{theorem} \label{fundthm}
Let $\cV$ be a vector space over a field $\F$. 
Suppose that $E$ is a linear endomorphism on
$\cV$, and $P=P[E]:\cV\to \cV$ is a linear operator polynomial in
$E$ which factors as
$$ P[E] = (E-\la_1) \cdots (E-\la_p) $$
where the scalars $\la_1,\ldots,\la_p \in \F$ are mutually distinct.
Then the solution space $\cV_P$, for $P$, admits a canonical
and unique direct sum decomposition
\begin{equation}
\cV_P=\oplus_{i=1}^\ell \cV_{\lambda_i}~,
\end{equation}
where, for each $i$ in the sum, $\cV_{\lambda_i}$ is the solution
space for $E-\lambda_i$. The projection $\Proj_i: \cV_P\to
\cV_{\lambda_i}$ is given by the formula
$$ \Proj_i = y_i \prod_{i \not= j=1}^{j=p} (E-\la_j)
   \quad \text{where} \quad
   y_i = \prod_{i \not= j=1}^{j=p} \frac{1}{\la_i - \la_j}. $$
\end{theorem} 

This is a special case of Theorem 1.1 from \cite{GoSiDec}.
To use this result in our setting we need the following result. 

\begin{lemma} \label{distinct}
The constants
$$
  - \frac{2i(n-2k-i+1)}{n}, \quad
  i \in \N,\ k \in \{0,\cdots ,n-1 \},
$$
 are mutually distinct and
negative for $i=1,\ldots,\frac{n-2k}{2}$.
\end{lemma}

\begin{proof}
Assume that the scalar 
$$ 2i(n-2k-i+1) - 2j(n-2k-j+1) = 2(i-j)[n-2k-(i+j)+1] $$
is equal to zero for some $i,j=1,\ldots,\frac{n-2k}{2}$.  This can
happen only if $i=j$ or $n-2k-(i+j)+1=0$. But the latter possibility
cannot happen as $i,j \leq \frac{n-2k}{2}$ means $i+j \leq n-2k$. Thus
the discussed scalars are mutually distinct. The scalars are negative
since for the ranges considered $i\geq 1$, and $2i\leq n-2k$ implies that
$i+2k\leq n-i<n+1$.
\end{proof}

\noindent{\bf Definition:} We define the scalars: 
\begin{equation} \label{lambda}
 \la^k_i:=  - \frac{2i(n-2k-i+1)}{n} J, \quad
  i \in \N,\ k \in \{0,\cdots ,n-1 \},
\end{equation}
where, recall, $J$ is the trace of the Schouten tensor. So on Einstein
manifolds these scalars are constant and, if $J\neq 0$ then these are
non-zero and mutually distinct.

\begin{proposition} \label{nullL} Let $(M,g)$ be an Einstein manifold
  which is not Ricci flat.  We will use the scalars $\la_i^k$ from
  \nn{lambda} and put $p=\frac{n-2k}{2}$.  The null space of the
  conformally invariant operator $L_k: \cE^k \to \cE_k$ defined in
  Corollary \ref{L_k} is
$$ \cN(L_k) = \widetilde{\cH}^k_{\si,0} \oplus \bigoplus_{i=1}^{p-1}
   \widetilde{\cH}^k_{\si,\la_i^{k+1}}, \quad k \in
   \{0,\ldots,\frac{n}{2}-1\}, $$ and $\cN(L_{n/2}) = \cE^{n/2} =
   \cE_{n/2}$.
\end{proposition}

\begin{proof}
The case $k=n/2$ is obvious so assume $k \in \{0,\ldots,\frac{n}{2}-1\}$.
Since $L_k = \de P^{p-1}_{k+1}[d\de] d$ according to Corollary 
\ref{L_k} and $P^{p-1}_{k+1}[d\de] d = d P^{p-1}_{k+1}[\de d]$ we get
$L_k = \de d P^{p-1}_{k+1}[\de d]$.
Now the Proposition follows from Theorem \ref{fundthm} for $E = \de d$ and 
Lemma \ref{distinct}.
\end{proof}
\noindent Note that $\cC^k\subseteq \widetilde{\cH}^k_{\si,0}$,
$\cC^k\cap \oplus_{i=1}^{p-1} \widetilde{\cH}^k_{\si,\la_i^{k+1}}=\{0\}$,
and $\oplus_{i=1}^{p-1} \widetilde{\cH}^k_{\si,\la_i^{k+1}} \subseteq
\cR(\delta)$ for $k \not= \frac{n}{2}$.

\begin{lemma} \label{nullG} Let $(M,g)$ be an Einstein manifold which
  is not Ricci flat.  In the Einstein scale $\si$, the null space of
  the operator $G_k^\si: \cE^k \to \cE_{k-1}$ given by \nn{Gk} is
\begin{equation} 
  \cN(G_k^\si) = \cN(\de) \oplus
   \bigoplus_{i=1}^p \overline{\cH}^k_{\si,\la_i^k},
\end{equation}
where the scalars $\la_i^k$ are from \nn{lambda} and $p=\frac{n-2k}{2}$.
\end{lemma}

\begin{proof}
Observe $G_k^\si = \de P^p_k [d\de] = P^p_k [\de d] \de$.
We will work in the Einstein scale $\si$ throughout the proof.
The case $k=\frac{n}{2}$ is obvious so we assume $k < \frac{n}{2}$.

Let us start the inclusion $\supseteq$. Clearly 
$\cN(\de) \subseteq \cN(G_k^\si)$.
Further suppose that $f \in \overline{\cH}^k_{\si,\la_j}$ for some
$j \in 1,\ldots,\frac{n-2k}{2}$. This means $d\de f = \la_j^k f$ 
using the definition of $\overline{\cH}^k_{\si,\la_j^k}$. 
The composition factors in (\ref{Q_k}), which yields the formula 
for $P^p_k [d\de]$, commute and one of these factors is
$d\de - \la_j^k$. Hence
$P^p_k [d\de] f=0$ which means $f \in \cN(G_k^\si)$.

Now we discusses the inclusion 
$\cN(G_k^\si) \subseteq \cN(\de) \oplus 
\bigoplus_{i=1}^{p} \overline{\cH}_{\si,\la_i^k}^k$.
From Lemma \ref{distinct} it follows that  the $\la_i^k$ for 
$i = 1,\ldots,p$ are mutually distinct.
First observe
$\cN(G_k^\si) \subseteq \cN(d\de P^p_k [d\de]: \cE^k \to \cE^k)$
since $G_k^\si = \de P^p_k [d\de]$.
It follows from Theorem \ref{fundthm} (where we put $E = d\de$) that 
$f \in \cN(d\de P^p_k [d\de])$ can be uniquely written
in the form
$f = \bar{f} + \sum_{i=1}^p f_i$ where 
$$ \bar{f} = \bar{y} \Bigl[ \prod_{j=1}^p (d\de - \la_j^k) \Bigr] f
   \quad \mbox{and} \quad
   f_i = y_i \, d\de  \Bigl[ \prod_{i \neq j=1}^{j=p} (d\de - \la_j^k) \Bigr]f, 
   \ \ i=1,\ldots,\ell $$
where $\bar{y}$ and $y_i$ are appropriate scalars,   
and this decomposition satisfies  
$d\de f_i = \la_i^k f_i$ for $i=1,\ldots,p$ and $d\de \bar{f}=0$.
Note the last display means $\bar{f} = \bar{y} P^p_k[d\de] f$.

Now consider this decomposition for $f \in \cN(G_k^\si)$. 
The condition $d\de f_i = \la_i^k f_i$ means 
$f_i \in \overline{\cH}_{\si,\la_i^k}^k$.
Further applying $\de$ to $\bar{f} = \bar{y} P^p_k[d\de] f$, we obtain
$\de \bar{f} = \bar{y} \de P^p_k[d\de] f = \bar{y} G_k^\si f =0$.
Hence $\bar{f} \in \cN(\de)$.
\end{proof}

\begin{theorem} Let $(M,g)$ be an Einstein manifold which is not Ricci
  flat. We use $\la_i^k$ to denote the scalars from \nn{lambda} and
  put $p=\frac{n-2k}{2}$.  The conformally invariant space $\cH_G^k =
  \cN(G_k: \cC^k \longrightarrow \cE_{k-1})$ is given by the direct sum
\begin{equation} \label{decomp}
  \cH_G^k = \cH_\si^k \oplus
  \bigoplus_{i=1}^p \overline{\cH}^k_{\si,\la_i^k}.
\end{equation}
\end{theorem}
\begin{proof}
By the definition, $\cH_G^k$ is equal to the intersection 
$\cN(G_k^\si) \cap \cN(d)$. Since 
$\cN(G_k^\si) = 
\cN(\de) \oplus \bigoplus_{i=1}^p \overline{\cH}^k_{\si,\la_i^k}$
according to Lemma \ref{nullG}, and 
$\overline{\cH}^k_{\si,\la_i^k} \subseteq \cR(d) \subseteq \cN(d)$, 
we obtain 
$\cH_G^k = 
(\cN(\de) \cap \cN(d)) \oplus \bigoplus_{i=1}^p 
\overline{\cH}^k_{\si,\la_i^k}$. 
\end{proof}

As discussed in the Introduction, a conformally invariant cohomology
space may be defined:
$$
H^k_L:=\cN(L_k)/\cR(d)~,
$$ where of course $d$ means $d:\ce^{k-1}\to \ce^k$, as is clear by
context.  From the definitions of the various spaces it follows
automatically that this fits into the complex \nn{coco} of
\cite{BrGodeRham}.
Here and below we put $H^{-1} :=0$ and $H^{-1}_L :=0$. 

Let us work on a non Ricci-flat Einstein manifold.
To discuss \nn{coco}
we note that we have mappings
\begin{equation}\label{bi}
 \overline{\cH}^k_{\si,\la} \stackrel{\de}{\longrightarrow}
     \widetilde{\cH}^{k-1}_{\si,\la} \quad \mbox{and} \quad
   \widetilde{\cH}^{k-1}_{\si,\la} \stackrel{d}{\longrightarrow}
     \overline{\cH}^k_{\si,\la}. 
\end{equation}
In the case $\la \not =0$ this is a bijective correspondence
as $f \stackrel{\de}{\mapsto} \de f \stackrel{d}{\mapsto} d\de f = \la f$ for
$f \in \overline{\cH}^k_{\si,\la}$ and similarly for the opposite direction.
In particular 
$$
d: \bigoplus_{i=1}^{\frac{n-2k}{2}}
\widetilde{\cH}^{k-1}_{\si,\la_i^{k}} \to \bigoplus_{i=1}^{\frac{n-2k}{2}}
\overline{\cH}^k_{\si,\la_i^k},
$$
is a bijection.  Using this it is easily verified that our results
above are consistent with \nn{coco} in the sense that from Proposition
\ref{nullL} and \nn{decomp} one verifies that \nn{coco} is exact at
$\cH^k_G$. Also since $\oplus_{i=1}^{\frac{n-2k}{2}}
\widetilde{\cH}^{k-1}_{\si,\la_i^{k}}$ intersects trivially with
$\cR(d)$ it follows that, in the Einstein scale $\si$,
\begin{equation}\label{Hkl}
H^{k-1}_L= \bigoplus_{i=1}^{\frac{n-2k}{2}}
\widetilde{\cH}^{k-1}_{\si,\la_i^{k}} \oplus 
\big( \widetilde{\cH}^{k-1}_{\si,0}/\cR(d) \big)~.
\end{equation}

The spaces $H^{*}_L$ also contribute to another exact complex on even
conformal manifolds, 
\begin{equation}\label{coco2}
0\to H^{k-1}\to H^{k-1}_L\stackrel{d}{\to}\mathcal{H}^k_{\L} \to H^k_L~, \quad \quad k\in \{0,1,\cdots ,n/2\}
\end{equation}
which is a tautological consequence of the definition:
$\mathcal{H}^k_{\L} =\cN(\L_k)$. In analogy with \nn{coco},
surjectivity of the last map is not known in general.  This motivates
studying $\cN(\L_k) $. From the results for the null spaces of $L_k$ and
$G_k$ we have the following.

\begin{corollary} \label{nullLG} The null space $\cN(\L_k) $, of
  $\L_k$, is conformally invariant. On  $(M,g)$,  an Einstein
  manifold which is not Ricci flat, $\cN(\L_k) $ is 
  given by the direct sum:
$$ \cN(\L_k) = \bigl( \cN(\de) \cap \cN(\de d) \bigr) \oplus
   \bigoplus_{i=1}^{p-1} \widetilde{\cH}^k_{\si,\la_i^{k+1}} \oplus   
   \bigoplus_{i=1}^p \overline{\cH}^k_{\si,\la_i^k}, 
   \quad 0 \leq k \leq \frac{n}{2}-1$$
where   the scalars $\la_i^k$ are given in  \nn{lambda} and 
$p=\frac{n-2k}{2}$. Further $\cN(\L_{n/2}) = \cN(\de)$.
\end{corollary}

\begin{proof}
The case $k=0$ follows from Proposition \ref{nullL} and the case
$k=\frac{n}{2}$ is obvious. 
Assume $1 \leq k \leq \frac{n}{2}-1$.
  The conformal invariance follows from the conformally invariance of
  $\L_k$. In the scale $\si$, we have $\cN(\L_k) = \cN(L_k) \cap
  \cN(G_k^\si)$ hence we need intersection of the direct sums in
 the displays of Proposition \ref{nullL} and Lemma \ref{nullG}.  Since
  $\overline{\cH}^k_{\si,\la_i^k} \subseteq \cN(\de d) =
  \widetilde{\cH}^k_{\si,0}$ (by the definition of these spaces) for $i =
  1,\ldots,p$, we obtain
$$ \cN(L_k) \cap \cN(G_k^\si) = \cN(\de) \cap \Bigl[
   \widetilde{\cH}^k_{\si,0} \oplus
   \bigoplus_{i=1}^{p-1} \widetilde{\cH}^k_{\si,\la_i^{k+1}} \Bigr]. $$   
Since similarly $\widetilde{\cH}^k_{\si,\la_i^{k+1}} \subseteq \cN(\de)$ for
$i=1,\ldots,p-1$, the statement follows.
\end{proof}
Returning to the sequence \nn{coco2}, exactness is easily verified using
Corollary \ref{nullLG} and the bijections \nn{bi}.
Note the direct summand $\cN(\de) \cap \cN(\de d)$ from the previous 
corollary satisfies
\begin{equation} \label{between} 
  \cN(\de) \cap \cN(d) \subseteq \cN(\de) \cap \cN(\de d) \subseteq
  \cN(d\de + \de d) ~.
\end{equation}

\begin{remark*}
  We can also study the space $\cN(Q_k^\si:\cC^k\to \ce_k)$. 
A direct
  application of Theorem \ref{fundthm} shows that 
$$ \cN(Q_k^\si:\cC^k\to \ce_k) =  \bigoplus_{i=1}^p \overline{\cH}^k_{\si,\la_i^k} \subseteq \cR(d) ,$$
on $(M,g)$  an Einstein manifold which is not Ricci
  flat and with notation as above.
\end{remark*}

 Let us, as usual,
assume that $(M,g)$ is Einstein and not Ricci flat. 
The operator $Q_k^\si$ simplifies on $\cH^k_G$. Considering \nn{decomp},
observe that
$Q_k^\si$ vanishes on $\overline{\cH}^k_{\si,\la^k_{i}} \subseteq
\cH^k_G$,  for each of the nonzero scalars $\la^k_i$;
this is because the composition factor $d\de - \la^k_{i}$ of $Q_k^\si$
vanishes on $\overline{\cH}^k_{\si,\la^k_{i}}$. Further, $Q_k^\si$ is
a multiple of the identity on $\cH^k_{\si}$ because $\de$ vanishes on
$\cH^k_{\si}$.  Using \nn{decomp}, we summarise this.
\begin{proposition}\label{Qres} 
Let $(M,g)$ be an Einstein manifold which is not Ricci flat.
The restriction of 
$Q_k^\si:\cC^k\to \ce_k$ to the
  conformal harmonics $\cH^k_G$ is given in the
  Einstein scale $\si$ as follows:
\begin{gather*}
  Q_k^\si|_{\overline{\cH}^k_{\si,\la^k_{i}}} =0 \quad \mbox{and} \quad
  Q_k^\si|_{\cH^k_{\si}} = s^k J^{(n-2k)/2}\id \ \ \mbox{where} \\
  s^k = \prod_{i=1}^{\frac{n-2k}{2}} \frac{2i(n-2k-i+1)}{n}.
\end{gather*}
\end{proposition}
\noindent
Note the last display means $s^{n/2}=1$.

The space $\cB^k = \{ df \mid Q_k^\si df \in \cR(\de) \} \subseteq
\cH^k_G$ is conformally invariant and in
\cite{BrGodeRham} plays a role in studying $\cH^k_G$. Clearly $f': =df
\in \overline{\cH}^k_{\si,\la_i^k}$, $i = 1,\ldots,p :=
\frac{n-2k}{2}$ satisfies $Q_k^\si f'=0$ thus trivially $f' \in
\cB^k$. Therefore $\bigoplus_{i=1}^p \overline{\cH}^k_{\si,\la_i^k}
\subseteq \cB^k.$

\section{Compact conformally-Einstein spaces}\label{compact}

Recall that for $\varphi\in {\mathcal{E}}^k$, $\psi\in
{\mathcal{E}}_k $, and $(M,[g])$ compact of signature $(p,q)$, there
is the natural conformally invariant global pairing
\begin{equation}\label{pair}
\varphi,\psi \mapsto \langle \varphi,\psi \rangle :=\int_M \varphi\cdot \psi\,
d\mu_{\mbox{\scriptsize\boldmath{$g$}}},
\end{equation}
where $\varphi\cdot \psi\in{\mathcal{E}}[-n]$
denotes a complete contraction between $\varphi $ and $\psi$. When $M$ is orientable we have
$$
\langle \varphi,\psi \rangle =\int_M \varphi \wedge \star \psi 
$$
where $\star $ is the conformal Hodge star operator.

This pairing combines with the operator $Q^\si_k$ to yield other
global pairings. For example on compact pseudo-Riemannian manifolds,
there is a conformally invariant pairing between
${\mathcal N}(L_k)$ and ${\mathcal{C}}^k$ given by
\begin{equation}\label{Qp}
(u,w)\mapsto \langle u, Q_k w \rangle \quad \quad k=0,1,\cdots ,n/2
\end{equation}
for $w\in{\mathcal{C}}^k $ and $u\in {\mathcal{N}}(L_k)$ \cite[Theorem
2.9,(ii)]{BrGodeRham}. We want to examine this in the Einstein
setting.  The case $k=n/2$ just recovers \nn{pair} and so we focus on
the remaining cases.

We assume that $(M,g)$ is even dimensional, Einstein and not
Ricci-flat. Consider $ \langle u, Q_k w \rangle $ with $w\in \cC^k$,
$u\in \cN(L_k)$, $k\in \{0,\cdots n/2-1\} $. By Proposition \ref{nullL}
$u$ decomposes directly:
$u=u_0+u_1$ where $u_0\in \widetilde{\cH}^k_{\si,0}$ and $u_1\in
\oplus_{i=1}^{p-1} \widetilde{\cH}^k_{\si,\la_i^{k+1}}.$
Now $u_1=\delta u'$ for some $(k+1)$-form $u'$ so, integrating by parts, 
$$
\langle u_1, Q^\si_k w \rangle = \langle u', d Q^\si_k w \rangle~.
$$ 
But using Corollary \ref{QCC} we have $d Q^\si_k w=0$. We summarise
 this simplification of the pairing.
\begin{lemma}\label{simpl}
  On an even, Einstein, and non Ricci-flat, compact manifold
  $(M,g=\si^{-2}\bg)$ the pairing on ${\mathcal N}(L_k)\times
  {\mathcal{C}}^k$ descends to $\widetilde{\cH}^k_{\si,0}\times
  \cC^k$.
\end{lemma}
\noindent Note that for $k=0$, $\widetilde{\cH}^k_{\si,0}$ is the 
null space of the Laplacian.

Next, on compact pseudo-Riemannian manifolds, we recall that
$Q_k^{\sigma}$ also gives a conformally invariant quadratic form
on $\cH^k_G$ \cite{BrGocomm}; this is given by \nn{Qp} with now $u,w\in
\cH^k_G$. We write this as
\begin{equation}\label{thetap}
\tilde{\Theta}:\cH^k_G\times \cH^k_G\to \mathbb{R}~.
\end{equation}
By Proposition \ref{Qres}, this specialises as follows:
\begin{proposition} \label{Qdes}
On non Ricci-flat  compact even Einstein manifolds $(M,g)$ the quadratic form 
$\tilde{\Theta}:\cH^k_G\times \cH^k_G\to \mathbb{R}$ descends to 
$$
\cH^k_\si\times \cH^k_\si \to \mathbb{R} \quad \quad k\in \{0,1,\dots ,n/2\}
$$
given by 
$$
(u,w)\mapsto s^k J^{\frac{n-2k}{2}}\langle u, w \rangle
$$
where the constant $s^k$ is given in Proposition \ref{Qres}.
\end{proposition}

\subsection{Compact Riemannian spaces}

We now assume $(M,g)$ is a compact Einstein manifold of Riemannian
signature. As above we relate $g$ to $\si\in \ce_+[1]$ by
$g=\si^{-2}\bg$,  we assume $k \in \{0,\ldots,\frac{n}{2}\}$
and we set $\cE_{-1}:=0$.  
Many results from the previous section simplify in
this setting.  In particular we may use the de Rham Hodge
decomposition $\cE^k = \cR(d) \oplus \cR(\de) \oplus \cH^k_\si$ and
$\cH^k_\si$ is the usual space of de Rham  harmonics, that is $\cH^k_\si =
\cN(d\de + \de d) \cong H^k$. It also follows that the containments in
\nn{between} may be replaced by set equalities.  Note also that, for
example, $\cN(\delta d)=\cN(d)$.

 Next observe that, since the operators $\delta d$ and $d \delta$ are
positive, we have the following from Lemma \ref{distinct}.
\begin{proposition}\label{pn}
If $(M,g)$ is a positive scalar curvature compact Riemannian Einstein manifold
then 
$$
\widetilde{\cH}^{k'}_{\si,\la_i^{k}} =0 \quad \mbox{and} \quad \overline{\cH}^{k'}_{\si,\la_i^{k}} =0 \quad\quad  k'\in \{ 0, \cdots ,n\}
 $$
 for the  $\la_i^k$ as in
  \nn{lambda}. 
\end{proposition}

Using \nn{between} and the related observations  we have the following specialisations
of the results of Section \ref{null}.

\begin{theorem}\label{crnrf} Let $(M,g)$ be a compact Riemannian Einstein
manifold of even dimension. 
We have the exact sequences  
$$
 0\to H^{k-1}\to H^{k-1}_L\stackrel{d}{\to}\mathcal{H}^k_G \to H^k\to 0,
$$
and 
$$
0\to H^{k-1}\to H^{k-1}_L\stackrel{d}{\to}\mathcal{H}^k_{\L} \to H^k_L \to 0.
$$
In particular $(M,[g]) $ is $(k-1)$-regular for $k=1,\cdots ,n/2$.

Assume $(M,g)$ is not Ricci-flat. With the scalars $\la_i^k$ as in
\nn{lambda} and $p=\frac{n-2k}{2}$ we have:
$$
 \cN(L_k) = \cC^k \oplus
   \bigoplus_{i=1}^{p-1} \widetilde{\cH}^k_{\si,\la_i^{k+1}},
   \quad k < n/2.
$$
$$
\cH_G^k = \cH_\si^k \oplus
  \bigoplus_{i=1}^p \overline{\cH}^k_{\si,\la_i^k}.
$$

$$ \cN(\L_k) = \cH^k_\si \oplus
   \bigoplus_{i=1}^{p-1} \widetilde{\cH}^k_{\si,\la_i^{k+1}} \oplus
   \bigoplus_{i=1}^p \overline{\cH}^k_{\si,\la_i^k}~,
   \quad k < n/2, 
$$
while trivially we have $\cN(L_{n/2}) = \cE^{n/2}$ and 
$\cN(\L_{n/2}) = \cN(\de)$.
In particular, in the case of positive scalar curvature:  
$\cN(L_k) = \cC^k$ and $\cN(\L_k) = \cH^k_\si$ for $k < n/2$, 
and $\cH_G^k = \cH_\si^k$.

If $(M,g)$ is Ricci-flat then
$$ \cN(L_k) = \cC^k \quad \text{and} \quad
   \cH_G^k = \cN(\L_k) = \cH^k_\si
   \quad \text{for} \quad k < n/2. $$
\end{theorem}

\noindent Note that in the non Ricci-flat case
$\cH_G^k$ is formally as in \nn{decomp}, but now we have $\cH_\si^k\cong H^k$. 

The implications for the global
pairings are as follows. The first statement in the Theorem follows from 
Lemma \ref{simpl}, Theorem \ref{crnrf},  and that 
$\widetilde{\cH}^k_{\si,0}=\cC^k$ in
the compact Riemannian setting.  
\begin{theorem} \label{Rpairs} Let $(M,g)$ be a compact Riemannian
   Einstein manifold of even dimension.  The pairing on ${\mathcal
     N}(L_k)\times {\mathcal{C}}^k$, by $ (u,w)\mapsto \langle u, Q_k
   w \rangle$, with $ k=0,1,\cdots ,n/2 $ descends to $\cC^k \times
   \cC^k$.  By Theorem \ref{crnrf}, the quadratic form $\tilde{\Theta}$
   from \nn{thetap} yields a conformally invariant
   quadratic form
$$
H^k\times H^k\to \mathbb{R}. 
$$
In the Einstein scale this is given by Proposition \ref{Qdes} where
$\cH^k_\si$ are the usual harmonics for $g$.
In the Ricci--flat case this quadratic form is zero for $k<n/2$ and 
recovers \nn{pair} for $k=n/2$.
\end{theorem}
\noindent The last statement of the Theorem uses expression \nn{Q_k}
and Theorem \ref{crnrf}. 
\begin{remark}
  Note that, for  the case of $k=0$ and $M$ connected, the first result of
  the Theorem states that for $f$ in the null space of the dimension
  order GJMS operator (recall
  $L_0=\Delta^{n/2}+lower~order~terms$)
$$
\int f Q =c\int Q.
$$
where $c$ is a unique constant such that $c-f\in \cR(\delta^\si)$.
(Here we write $\delta^\si$ to emphasise that, although the display is
conformally invariant, to write the difference $c-f$ as a divergence
requires working in the Einstein scale.)
\end{remark}

\medskip

\begin{corollary}
Put $p := \frac{n-2k}{2}$. In an Einstein scale the space $\cB^k$ is given as follows:
$$ \cB^k =
   \begin{cases} \oplus_{j=1}^{p}
                 \overline{\cH}^k_{\si,\la_i^k} & J \not= 0 \\
                 0 & J=0.
   \end{cases} $$
\end{corollary}

\section{The Fefferman-Graham ambient metric}\label{asec}

Thus let us review briefly the basic relationship between the
Fefferman-Graham ambient metric construction and tractor calculus as
described in \cite{CapGoamb} for general conformal manifolds.

Let $\pi:\cq\to M$ be a conformal structure of signature $(p,q)$.
Let us use $\rho $ to denote the ${\Bbb R}_+$ action on $ \cq$ given
by $\rho(s) (x,g_x)=(x,s^2g_x)$.  An {\em ambient manifold\/} is a
smooth $(n+2)$-manifold $\aM$ endowed with a free $\Bbb R_+$--action
$\rho$ and an $\Bbb R_+$--equivariant embedding
$i:\cq\to\aM$.  We write $\aX\in\frak X(\aM)$ for the fundamental field
generating the $\Bbb R_+$--action, that is for $f\in C^\infty(\aM)$
and $ u\in \aM$ we have $\aX f(u)=(d/dt)f(\rho(e^t)u)|_{t=0}$.

If $i:\cq\to\aM$ is an ambient manifold, then an {\em ambient
metric\/} is a pseudo--Riemannian metric $\h$ of signature $(p+1,q+1)$
on $\aM$ such that the following conditions hold:

\smallskip
\noindent
(i) The metric $\h$ is homogeneous of degree 2 with respect to the
$\Bbb R_+$--action, i.e.\ if $\Cal L_{\sX}$
denotes the Lie derivative by $\aX$, then we have $\Cal L_{\sX}\h=2\h$.
(I.e.\ $\aX$ is a homothetic vector field for $h$.)

\noindent
(ii) For $u=(x,g_x)\in \cq$ and $\xi,\eta\in T_u\cq$, we have
$\h(i_*\xi,i_*\eta)=g_x(\pi_*\xi,\pi_*\eta)$.

\noindent
To simplify the notation  we will usually identify $\cq$
with its image in $\aM$ and suppress
the embedding map $i$.

\smallskip

To link the geometry of the ambient manifold to the underlying
conformal structure on $M $ one requires further conditions. In
\cite{FGast,FGrNew} Fefferman and Graham treat the 
construction of a formal power series solution, along $ \Cal Q$, for the
Goursat problem of finding an ambient metric $ \h$ satisfying (i) and
(ii) and the condition that it be Ricci flat, i.e.\ Ric$(\h)=0$. In
even dimensions for a general conformal structure this is obstructed
at finite order.  However when the underlying conformal structure is
(conformally) Einstein then an explicit Ricci-flat ambient metric is
available \cite{GrH,Leis,Leit}.  (In fact also more generally a
similar result is available for certain products of Einstein manifolds
\cite{GoL}.)  Here we shall use only the existence part of Ricci-flat
ambient metric. The uniqueness of the operators we will construct is a
consequence of the fact that they can be uniquely expressed in terms
of the underlying conformal structure as in \cite{CapGoamb,GoPetLap}.

It turns out that one may arrange that $\h$ is a metric satisfying the
conditions above (i.e.\ (i) and (ii) and with $\h$ Ricci flat to the
order possible) with $ Q:=\h(\aX,\aX)$  a defining
function for $\cq$, and $2\h(\aX,\cdot)=dQ$ to all orders in both odd
and even dimensions.  We write $ \nda $ for the ambient Levi-Civita
connection determined by $ \h$.  We will use upper case abstract
indices $A,B,\cdots $ for tensors on $ \aM$. For example, if $ v^B$ is
a vector field on $\aM $, then the ambient Riemann tensor will be
denoted $\aR_{AB}{}^C{}_{D}$ and defined by $
[\nda_A,\nda_B]v^C=\aR_{AB}{}^C{}_{D}v^D$. In this notation the
ambient metric is denoted $ \h_{AB}$ and, with its inverse, this is
used to raise and lower indices in the usual way. We will not normally
distinguish tensors related in this way even in index free notation;
the meaning should be clear from the context.  Thus for example we
shall use $\aX$ to mean both the Euler vector field $\aX^A$ and the
1-form $ \aX_A=\h_{AB}\aX^B$.

Let $ \cce(w)$ denote the space of functions on $ \aM$ which are
homogeneous of degree $ w\in {\Bbb R}$ with respect to the action $
\rho$. That is $f\in \cce(w)$ means that $ \aX f=wf$. Similarly a
tensor field $F$ on $ \aM$ is said to be {\em homogeneous of degree}
$w$ if $\rho(s)^* F= s^w F$ or equivalently $ \cL_{\miniX} F=w F$.
Just as sections of $ \ce[w]$ are equivalent to functions in $
\cce(w)|_\cq$ we will see that the restriction of homogeneous tensor
fields to $\cq$ have interpretations on $M$ as weighted sections of 
tractor bundles \cite{CapGoamb,GoPetLap}.

On the ambient tangent bundle $T\aM$ we define an action of $\Bbb R_+$
by $s\cdot \xi:=s^{-1}\rho(s)_\ast \xi$.  The sections of $ T\aM$
which are fixed by this action are those which are homogeneous of
degree $ -1$. Let us denote by $ \act$ the space of such sections and
write $\act(w)$ for sections in $\act\otimes \cce(w)$, where the
$\otimes$ here indicates a tensor product over $\cce(0)$. 
Along $ \cq$ the $\Bbb
R_+$ action on $T\aM$ is compatible with the ${\Bbb R}_+$ action on
$\cq$, so the quotient $(T\aM|_\cq)/\Bbb R_+$, is a rank $ n+2$ vector
bundle over $\cq/\Bbb R_+=M$; in fact this is (up to isomorphism) the
normal standard tractor bundle $\cT$ (or $\ce^A$) \cite{CapGoamb,GoPetLap} and the
composition structure of $\cT$ reflects the vertical subbundle 
$T\cq$ in $T\aM|_\cq$.  Sections of $\cT$ are equivalent to sections
of $T\aM|_\cq$ which are homogeneous of degree $-1$, that is sections of 
$\act|_\cq$. Using this
relationship one sees that the ambient metric $\h$ and the ambient
connection $\nda$ descend to, respectively the tractor metric $h$, and
the tractor connection $\nd^\cT$. For the metric this is obvious. We
discuss the connection briefly.  For $ U\in \ct$, let $ \tU $ be the
corresponding section of $\act|_\cq$.
A tangent vector field $ \xi$ on $M$ has a lift to a
homogeneous degree 0 section $ \tilde{\xi}$, of $T\aM|_\cq$, which is
everywhere tangent to ${\Cal Q}$.  This is unique up to adding $f\aX$,
where $ f\in \cce(0)|_\cq$. We extend $ \tU$ and $\tilde{\xi} $
smoothly and homogeneously to fields on $ \aM$.  Then we can form $
\nda_{\tilde{\xi}} \tU$; along $\cq$, this is clearly independent of
the extensions. Since $ \NX \tU =0$, the section $ \nda_{\tilde{\xi}}
\tU$ is also independent of the choice of $ \tilde{\xi}$ as a lift of
$ \xi$. Finally, the restriction of $ \nda_{\tilde{\xi}} \tU$ is a
homogeneous degree $-1$ section of $T\aM|_\cq $ and so determines a
section of $ \cT$ which depends only on $ U$ and $ \xi$.  This is
$\nd^\cT U$.

Finally we will say that an ambient tensor $F$ is homogeneous of {\em
  weight} $w$ if $\nda_{\sX}F=w F$. The weight is a convenient
shifting of homogeneity degree. Note, for example, that an ambient
1-form $\tilde{U}$ which is homogeneous of degree $-1$ is homogeneous
of weight 0 and this means that $\nda_{\sX}\tilde{U}=0 $.

\subsection{The main result}

In Section \ref{detour} 
several operators were defined on
conformally Einstein manifolds directly using tractor calculus and the
parallel tractor of the Einstein structure. On the other hand in
\cite{BrGodeRham} operators with the same notation were defined on
general conformal manifolds via the Fefferman-Graham ambient metric,
and its link to tractor calculus. The aim of this section is simply to
show that these agree (up to a nonzero multiple).

\begin{proposition} \label{ambLLk}
Assume $n$ even and $k \in \{ 1,\ldots,\frac{n}{2} \}$.
On Einstein manifolds the operator $\L_k$ defined by \nn{LLk}
 agrees  with the operator with the same
notation in \cite{BrGodeRham}. The operators $L_k$ and $G_k^\si$ from
Section \ref{detour} also agree with the operators of the same notation in
\cite{BrGodeRham}.
\end{proposition}

\noindent Here ``agree'' means the operator is the same up to a
non-zero multiple, and we will not pay attention to the detail of what
this constant factor is.

On the ambient manifold a special role is played by
differential operators $P$ on ambient tensor bundles which act {\em
  tangentially} along $\cq$, in the sense that $PQ=QP'$ for some
operator $P'$ (or equivalently $[P,Q]=QP''$ for some $P''$). Note that
compositions of tangential operators are tangential. If tangential
operators are homogeneous (i.e.\ the commutator with the Lie derivative
$\cL_{\sX}$ recovers a constant multiple of the operator) then they
descend to operators on $ M$. An example of a tangential operator is
given by
$$
(n+2w-2)\nda+\aX \afl=:\afD: \act^\Phi(w)\to \act\otimes  \act^\Phi(w-1)
$$
where $\act^\Phi(w)$ indicates the  space of sections,
homogeneous of weight
$w$,  of some ambient tensor bundle,  and  
$$
\afl={\mbox{\boldmath$\Delta$}}-\aR\sharp\sharp .
$$ 
Here we use the Laplacian ${\mbox{\boldmath$\Delta$}} := -\nda^A \nda_A$
for compatibility with \cite{BrGodeRham}.
We will leave the verification that $\afD$ is tangential to the
reader, but note this also follows from the result that
$(n+2w-2)\nda+\aX \mbox{\boldmath$\Delta$} = :\aD: \act^\Phi(w)\to
\act\otimes \act^\Phi(w-1)$ is tangential as discussed in
\cite{GoPetLap,CapGoamb}.  Since this is tangential and homogeneous it
descends to an operator on weighted tractors. In fact it gives the
usual tractor-D operator \cite{GoPetLap,CapGoamb}. The ambient
$\aR\sharp\sharp$ similarly descends (in dimensions $n \neq 4$) to a
multiple of $W\sharp\sharp$. Thus acting on weighted tractor bundles
\cite{GoPetobst}.  Thus $ \act^\Phi(w)\to \act\otimes
\act^\Phi(w-1)$ descends to $\ct^\Phi(w)\to \ct\otimes \ct^\Phi(w-1)$
in dimensions other than 4. (Here $\ct^\Phi$ means the tractor bundle corresponding to $\act^\Phi$.) Henceforth for $(M,[g])$ of dimension 4
we take $\afD:=\aD$, rather then the definition above.

Now if $(M,[g])$ is conformally Einstein and $I$ a parallel tractor
corresponding to an Einstein scale then along $\cq$ in $\aM$ we have a
corresponding parallel vector field $\aI$. From the explicit formula
for the ambient metric over an Einstein manifold ones sees that
$\aI$ extends to a parallel vector field on $\aM$. (In fact when the
Einstein scale is not Ricci flat then the ambient metric is given as a
product of the metric cone with a line.)  
We have (on $\act^\Phi[w]$)
$$
\aI^A\afD_A=(n+2w-2) \aI^A \nda_A +\bsi \afl~,
$$
where $\bsi=\aI_A \aX^A\in \cce(1)$. Note that $\bsi$ is a
homogeneous function on $\cq$ corresponding to $\si=I_A X^A$.

Thus if we extend a tensor field $\aU\in \act^\Phi(w)|_\cq$ off $\cq$
in  such a way that $\aI^A\nda_A \aU=0$ (which implies $\aU\in
\act^\Phi(w)$) then we get simply
$$
\aI^A\afD_A= \bsi \afl~.
$$
Note that $\aI^A\nda_A \aU=0$ can be achieved by starting with a section
along $\cq$ and then extending off $\cq$ by parallel transport. The
key point here is that $I^AX_A$ is non-vanishing, at least in a
neighbourhood of $\cq$, and so $\aI^A \nda_A$ is not tangential to
$\cq$.

Next observe that, since $\bsi=\aI_A \aX^A$ and $\aI_A$ is parallel, we have
$$
\nda_A \bsi = \aI_A~,
$$
which is parallel.
Thus
\begin{equation}\label{scalelap}
[\afl,\bsi]=[\aDelta,\bsi]=2\aI^A\nda_A
\end{equation}
where we consider $\bsi$ as a multiplication operator.

The following observations will be useful.
\begin{lemma}\label{IR0}
If $\aR$ denotes the ambient curvature then
 $\aI^A\nda_A \aR=0$.
\end{lemma}
\begin{proof}
By the Bianchi identity
$$
\aI^A\nda_A \aR_{BC}{}^D{}_E+\aI^A\nda_C \aR_{AB}{}^D{}_E+\aI^A\nda_B \aR_{CA}{}^D{}_E=0.
$$
But $\aI$ is parallel which implies that $[\nda,\aI]=0$ and
$\aI^A\aR_{AB}{}^D{}_E=0= \aI^A\aR_{CA}{}^D{}_E$.  So the result follows.
\end{proof}

\begin{lemma}\label{ndIafl}
If $\aU$ is an ambient tensor such that $\aI^A\nda_A \aU=0$ then, 
for any $p\in \N\cup \{0\}$,
$\aI^A\nda_A (\afl^p \aU)=0$
\end{lemma}
\begin{proof}
  Clearly, acting on any ambient tensor, we have
  $[\aI^A\nda_A,\nda_B]=0$.  Thus $\aI^A\nda_A$ commutes with the
  Bochner Laplacian $\aDelta$.  On the other hand by definition $\afl$ differs from  the Bochner
by a  curvature action:
$
\afl-{\mbox{\boldmath$\Delta$}}=-\aR\sharp\sharp ,
$
while from the previous Lemma the ambient curvature is parallel along
the flow of $\aI^A\nda_A$.
\end{proof}

The main technical result we need is this.
\begin{proposition}\label{comp}
For $\af$ an ambient form homogeneous of weight $k-n/2$ we have
$$
(\aI^A\afD_A)^k \af =\bsi^k \afl^k \af ~,
$$
along $\cq$.
\end{proposition}
\begin{proof}
  First note that both sides are tangential operators. For the
  right-hand-side this is proved in \cite{BrGodeRham}. For the left-hand-side
  it holds simply because $\afD$ is tangential and $\aI$ is parallel
  on the ambient manifold.   So neither side can depend
on the transverse (to $\cq$) derivatives of the homogeneous $\af$.

Now the result is true if $k=1$. Also, calculating along $\cq$,
$$
(\aI^A\afD_A)^k \af=  (\aI^B\afD_B)^{k-1}\aI^A\afD_A \af
$$
and so by induction
$$
(\aI^A\afD_A)^k \af=\bsi^{k-1}\afl^{k-1}\aI^A\afD_A \af~.
$$
Since the result is independent of transverse derivatives we may choose the extension off $\cq$ to suit. Thus we assume without loss of generality that
$\aI^A \nda_A \af=0$. Then $\aI^A\afD_A \af=\bsi \afl \af$
and so
$$
\bsi^{k-1}\afl^{k-1}(\aI^A\afD_A) \af=\bsi^{k-1}\afl^{k-1}(\bsi \afl \af) .
$$
So from \nn{scalelap} and Lemma \ref{ndIafl} the result follows.
\end{proof}

By Proposition 3.2 of \cite{BrGodeRham}, the operator $\afl^m $ is
homogeneous and acts tangentially  on ambient
differential forms of weight $m-n/2$. Thus it descends to an operator
that we denote $\fl_m$ on form-tractors of weight $m-n/2$. From the
above Proposition we obtain immediately the following results.
\begin{corollary}\label{flp}
  On conformally Einstein manifolds $(M,[g])$ the invariant operator
  $\fl_m:\cT^k[m-n/2]\to \cT^k[-m-n/2]$, $m\in\{0,1,2,\cdots \}$, is
  formally self-adjoint and given by
$$
\fl_m=
\si^{-m} (I^A\fD_A)^m
$$
where $\si^{-2}\bg$ is an Einstein metric on $M$ and
$I=\frac{1}{n}D\si$.
In odd dimensions these are natural operators.
In even dimensions the same is true with the restrictions that
either
$m\leq n/2-2$; or
$m\leq n/2-1$ and $k=1$; or
$m\leq n/2 $  and $ k=0 $.
In the conformally flat case the operators are natural with no
restrictions on $m\in \{1,2,\ldots\}$.
\end{corollary}

\begin{proof}
The statements on naturality are extracted from \cite{BrGodeRham}.
It only remains to establish the claim that the operator is formally
self-adjoint. But this is immediate from the formula for the right-hand-side 
from \nn{EmodBox} because $I^A \modD_A = \modBox_\si$ according to 
\nn{modBox}.
\end{proof}

Finally we are ready to prove the main result: \\
\newcommand{\bL}{\mathbb{L}}
\let\e=\varepsilon
\let\i=\iota
\noindent {\it Proof of Proposition \ref{ambLLk} :}
By expression (40) from \cite{BrGodeRham} and the fact that $\fl_m$, as in Corollary \ref{flp}, is formally self-adjoint we have
that the operator $\bL_k$ from \cite{BrGodeRham} is given by
$$\bL_k:= \fl_\ell\i(\fD)\e(X) q_k~,
$$ where the notation is from that source.  But it is a straightforward 
calculation to verify that, up to a non-zero
multiple, $\i(\fD)\e(X) q_k$ is exactly the operator $M$ from \nn{M}.
(See also \cite[2.1.2 and (2.8)]{Sthesis} where the special case $k=2$ is
treated in detail.)  
 So the result now follows from the Corollary and \nn{L}
where $w=0$ and $p = \frac{n-2k}{2}$.  \quad $\Box$

\end{document}